\documentclass[11pt]{amsart}
\usepackage[numbers, square]{natbib}
\usepackage{amssymb}
\usepackage{amsmath}
\usepackage{amsfonts}
\usepackage{geometry}
\usepackage{graphicx}
\usepackage{amsthm}
\usepackage{hyperref}
\usepackage[dvipsnames]{xcolor}

\providecommand{\U}[1]{\protect\rule{.1in}{.1in}}
\theoremstyle{plain}

\newtheorem{corollary}{Corollary}

\newtheorem{lemma}{Lemma}

\newtheorem{remark}{Remark}

\newtheorem{theorem}{Theorem}
\newtheorem*{theorem*}{Theorem}
\numberwithin{equation}{section}

\newcommand{\disp}{\displaystyle}

\newcommand{\eps}{\varepsilon}
\newcommand{\vp}{\varphi}

\newcommand{\al}{\alpha}
\newcommand{\be}{\beta}

\newcommand{\de}{\delta}

\newcommand{\te}{\theta}
\newcommand{\la}{\lambda}

\newcommand{\om}{\omega}
\newcommand{\Om}{\Omega}
\newcommand{\si}{\sigma}

\newcommand{\nid}{\noindent}

\newcommand{\iny}{\infty}

\newcommand{\su}{\subset}
\newcommand{\LP}{\Delta}
\newcommand{\gr}{\nabla}

\newcommand{\norm}[1]{\left\vert\left\vert #1\right\vert\right\vert}
\newcommand{\innp}[1]{\left< #1 \right>}
\newcommand{\abs}[1]{\left\vert#1\right\vert}
\newcommand{\set}[1]{\left\{#1\right\}}
\newcommand{\brac}[1]{\left[#1\right]}
\newcommand{\pr}[1]{\left( #1 \right) }
\newcommand{\pb}[1]{\left( #1 \right] }

\newcommand{\N}{\ensuremath{\mathbb{N}}}

\newcommand{\R}{\ensuremath{\mathbb{R}}}
\newcommand{\Z}{\ensuremath{\mathbb{Z}}}
\newcommand{\C}{\ensuremath{\mathbb{C}}}

\begin{document}
\title[Quantitative uniqueness ]
{Quantitative uniqueness of solutions to second order elliptic equations with
singular potentials in two dimensions}
\author{ Blair Davey and Jiuyi Zhu}
\address{
Department of Mathematics\\
The City College of New York\\
New York, NY 10031, USA\\
Email: bdavey@ccny.cuny.edu }
\address{
Department of Mathematics\\
Louisiana State University\\
Baton Rouge, LA 70803, USA\\
Email:  zhu@math.lsu.edu }
\thanks{\noindent{Davey is supported in part by the Simons Foundation Grant number 430198.
Zhu is supported in part by \indent NSF grant DMS-1656845}}
\date{}

\begin{abstract}
In this article, we study the vanishing order of solutions to second order elliptic equations with singular lower order terms in the plane. 
In particular, we derive lower bounds for solutions on arbitrarily small balls in terms of the  Lebesgue norms of the lower order terms for all admissible exponents.
Then we show that a scaling argument allows us to pass from these vanishing order estimates to estimates for the rate of decay of solutions at infinity. 
Our proofs rely on a new $L^p - L^q$ Carleman estimate for the Laplacian in $\mathbb{R}^2$.
\end{abstract}

\maketitle
\section{Introduction}

In this paper, we study the vanishing order of solutions to second order elliptic equations with singular lower order terms in $\R^2$.
Vanishing order is a characterization of the rate at which a solution can go to zero and is therefore a quantitative description of the strong unique continuation property.
Let $\Om \su \R^2$ be open and connected.
Assume that $V : \Om \to \C$ is an element of $L^t\pr{\Om}$ for some $t > 1$, and that $W : \Om \to \C^2$ belongs to $L^s\pr{\Om}$ for some $s > 2$.
Let $L = \LP + W \cdot \gr + V$ denote an elliptic operator.
Suppose $u$ is a solution to $L u = 0$ in $\Om$.
Our vanishing order estimates take the following form:
For all $r$ sufficiently small and $x_0 \in \Om' \Subset \Om$,
\begin{equation}
\norm{u}_{L^\iny\pr{B_r\pr{x_0}}} \ge c r^{C \be},
\label{vanishing}
\end{equation}
where $\be$ depends on the elliptic operator $L$ or, more explicitly, on the Lebesgue norms of $W$ and $V$.
 The estimate \eqref{vanishing} implies that the vanishing order of  $u$ at $x_0$ is less than $C\beta$.
Vanishing order estimates have a deep connection with the eigenfunctions Riemannian manifolds.
Let $\mathcal{M}$ be a smooth, compact  Riemannian manifolds without boundary, and let $ \phi_\lambda$ be a classical eigenfunction,
 $$-\LP_{g} \phi_\lambda=\lambda \phi_\lambda \quad \quad \mbox{in} \ \mathcal{M}.$$
Donnelly and Fefferman in \cite{DF88} proved that for any $x_0\in \mathcal{M}$
$$\norm{u}_{L^\iny\pr{B_r\pr{x_0}}} \ge c r^{C \sqrt{\lambda}}.$$
That is, the maximal vanishing order for eigenfunction $\phi_\lambda$ is $C\sqrt{\lambda}$ and this estimate is sharp.

More recently, there has been an interest in understanding how the vanishing order of solutions to the Schr\"odinger equation
\begin{equation}
-\LP u=V(x)u
\label{schro}
\end{equation}
depends on the potential function, $V$.
Bourgain and Kenig \cite{BK05} showed that when $V \in L^\infty$ the vanishing order is
less than $C(1+\|V\|_{L^\infty}^{2/3})$, i.e. the estimate \eqref{vanishing} holds with $\beta=(1+\|V\|_{L^\infty}^{2/3})$.
The examples constructed by Meshkov in \cite{Mes92} indicate that Bourgain and Kenig's result is sharp when $V$ is complex-valued and the underlying space is $2$-dimensional.
If $V \in W^{1, \infty}$, Kukavica in \cite{Kuk98} established that the upper bound for vanishing order is $C(1+\|V\|_{W^{1, \infty}})$.
By different approaches, Bakri \cite{Bak12} and Zhu \cite{Zhu16} independently improved Kukavica's results and showed that the sharp vanishing order is less than $C(1+\sqrt{\|V\|_{W^{1, \infty}}})$.

Since elliptic equations with singular lower order terms are known to satisfy the strong unique continuation property, attention has recently shifted to trying to understand the role of singular weights in vanishing order estimates.
Kenig and Wang \cite{KW15} characterized the vanishing order of solutions to elliptic equations with drift, $\LP u+W\cdot \nabla u = 0$, in the plane for real valued $W\in L^s(\mathbb R^2)$ with
$s\geq 2$.
Since we would like to consider possibly complex-valued $W$ and $V$, the techniques from \cite{KW15} are not applicable to the present paper where we investigate the vanishing order of solutions to $L u = 0$.
In \cite{DZ17}, the authors developed a new $L^p\to L^q$ type quantitative Carleman estimates for a range of $p$ and $q$ values.
This Carleman estimates allowed us to study the quantitative uniqueness  of solutions to $Lu=0$ for $n\geq 3$.
Motivated by the ideas in \cite{DZ17}, here we further develop the techniques to study solutions $Lu=0$ for $n=2$.

Through the application of a scaling argument, the vanishing order estimate given in \eqref{vanishing} implies a quantitative unique continuation at infinity theorem.
For solutions to $L u = 0$ in $\R^2$, the following result can be invoked from \eqref{vanishing}: For all $R$
sufficiently large,
\begin{equation}
\inf_{\abs{x_0} = R} \norm{u}_{L^\iny\pr{B_1\pr{x_0}}} \ge c \exp\pr{- C R^\Pi \log R},
\end{equation}
where $\Pi$ depends on $L$.
We may interpret such estimates as bounds on the rate of decay of solutions at infinity.
In \cite{BK05}, Bourgain and Kenig showed that if $u$ is a bounded, normalized, non-trivial solution to \eqref{schro} in $\R^n$, where $\norm{V}_{L^\iny\pr{\R^n}} \le 1$, then a vanishing order estimate implies that
\begin{align}
\inf_{\abs{x_0} = R} \norm{u}_{L^\iny\pr{B_1\pr{x_0}}} \ge c \exp\pr{- C R^{4/3} \log R}.\label{quantEst}
\end{align}
This lower bound played an important role in their study of Anderson localization.

Recall that $W \in L^s_{loc}\pr{\Om}$ and $V \in L^t_{loc}\pr{\Om}$ for some $s > 2$ and $t >1$.
Suppose $u\in W^{1,2}_{loc}\pr{\Om}$ is a non-trivial weak solution to
\begin{equation}
\LP u+ W(x)\cdot \nabla u+V(x) u=0.
\label{goal}
\end{equation}
An application of H\"older's inequality implies that there exists a $p \in \pb{1, 2}$, depending on $s$ and $t$, such that $W(x)\cdot \nabla u+V(x) u \in L^p_{loc}\pr{\Om}$.
By regularity theory, it follows that $u \in W_{loc}^{2,p}\pr{\Om}$ and therefore $u$ is a solution to \eqref{goal} almost everywhere in $\Om$.
Moreover, we have that $u \in L^\iny_{loc}\pr{\Om}$.
Therefore, the solution $u$ belongs to $L^\iny_{loc}\pr{\Om} \cap W^{1,2}_{loc}\pr{\Om} \cap
W^{2,p}_{loc}\pr{\Om}$ and $u$ satisfies equation \eqref{goal} almost everywhere in $\Om$.
We use the notation $B_{r}\pr{x_0} \su \R^n$ to denote the ball of radius $r$ centered at $x_0$. When the center is understood from the context, we simply write $B_{r}$.

Our first result, an order of vanishing estimate for equation \eqref{goal}, is as follows.

\begin{theorem}
Assume that for some $s \in \pb{2, \iny}$ and $t \in \pb{1, \iny}$, $\norm{W}_{L^s\pr{B_{10}}} \le K$ and $\norm{V}_{L^t\pr{B_{10}}} \le M$.
Let $u$ be a solution to \eqref{goal} in $B_{10}$.
Assume that $u$ is bounded and normalized in the sense that
\begin{align}
& \|u\|_{L^\infty(B_{6})}\leq \hat{C},
\label{bound} \\
& \|u\|_{L^\infty(B_{1})}\geq 1.
\label{normal}
\end{align}
Then for any sufficiently small $\eps > 0$, the vanishing order of $u$ in $B_{1}$ is less than $C\pr{1 + C_1 K^\kappa + C_2 M^\mu}$.
That is, for any $x_0\in B_1$ and every $r$ sufficiently small,
\begin{align*}
\|u\|_{L^\iny(B_{r}(x_0))}
 &\ge c r^{C\pr{1 + C_1 K^\kappa + C_2 M^\mu}},
\end{align*}
where
$\disp \kappa = \left\{\begin{array}{ll}
\frac{2s}{s-2} & t > \frac{2s}{s+2} \medskip \\
\frac{t}{t-1-\eps t} & 1 < t \le \frac{2s}{s+2}
\end{array}\right.$,
$\disp \mu = \left\{\begin{array}{ll}
\frac{2s}{3s-2} & t \ge s \medskip \\
\frac{2st}{3st+2t-4s-\eps(2st+4t-4s)}  & \frac{2s}{s+2} < t < s \medskip \\
\frac{t}{t-1+\eps\pr{t-2t\eps}} & 1 < t \le \frac{2s}{s+2}
\end{array}\right.$, \newline
$c = c\pr{s, t, K, M, \hat C, \eps}$, $C$ is universal, $C_1 = C_1\pr{s, t, \eps}$, and $C_2 = C_2\pr{ s, t,\eps}$.
\label{thh}
\end{theorem}

Theorem \ref{thh} in combination with a scaling argument gives rise to the following unique continuation at infinity theorem.

\begin{theorem}
Assume that for some $s \in \pb{2, \iny}$ and $t \in \pb{ 1, \iny}$, $\norm{W}_{L^s\pr{\R^2}} \le A_1$ and $\norm{V}_{L^t\pr{\R^2}} \le A_0$.
Let $u$ be a solution to \eqref{goal} in $\R^2$.
Assume that $\norm{u}_{L^\iny\pr{\R^2}} \le C_0$ and $\abs{u\pr{0}} \ge 1$.
Then for every sufficiently small $\eps > 0$ and all $R$ sufficiently large,
\begin{equation*}
\inf_{\abs{x_0} = R} \norm{u}_{L^\iny\pr{B_1\pr{x_0}}} \ge \exp\pr{-C R^\Pi \log R},
\end{equation*}
where
\begin{align*}
\Pi = \left\{\begin{array}{ll}
2 & t > \frac{2s}{s+2} \medskip \\
\frac{t(s-2)}{s(t-1-\eps t)} & 1 < t \le \frac{2s}{s+2},
\end{array}\right.
\end{align*}
and $C = C\pr{s, t, A_1, A_0, C_0, \eps}$.
\label{UCVW}
\end{theorem}

If we consider elliptic equations with drift, then the order of vanishing estimates follow directly from Theorem \ref{thh}.
In particular, if $V \equiv 0$, then the following statement is a consequence of Theorem \ref{thh} with $t = \iny$ and $M = 0$.

\begin{corollary}
Assume that for some $s \in \pb{2, \iny}$, $\norm{W}_{L^s\pr{B_{10}}} \le K$.
Let $u$ be a solution to $\LP u + W \cdot \gr u = 0$ in $B_{10}$.
Assume that $u$ is bounded and normalized in the sense of \eqref{bound} and \eqref{normal}. Then the vanishing order of $u$ in $B_{1}$ is less than $C\pr{1 + C_1 K^\kappa}$.
That is, for any $x_0\in B_1$ and every $r$ sufficiently small,
\begin{align*}
\|u\|_{L^\iny(B_{r}(x_0))}
 &\ge c r^{C\pr{1 + C_1 K^\kappa}},
\end{align*}
where
$\disp \kappa = \frac{2s}{s-2} $,
$c = c\pr{s, K, \hat C}$, and $C_1 = C_1\pr{s}$.
\label{thhW}
\end{corollary}

The corresponding unique continuation at infinity theorem follows from Theorem \ref{UCVW} in the same way.
Note that this pair of corollaries is independent of $\eps$.

\begin{corollary}
Assume that for some $s \in \pb{2, \iny}$, $\norm{W}_{L^s\pr{\R^n}} \le A_1$.
Let $u$ be a solution to $\LP u + W \cdot \gr u = 0$ in $\R^2$.
Assume that $\norm{u}_{L^\iny\pr{\R^n}} \le C_0$ and $\abs{u\pr{0}} \ge 1$.
Then for every $R$ sufficiently large,
\begin{equation*}
\inf_{\abs{x_0} = R} \norm{u}_{L^\iny\pr{B_1\pr{x_0}}} \ge \exp\pr{-C R^2 \log R},
\end{equation*}
where $C = C\pr{s, A_1, C_0}$.
\label{UCW}
\end{corollary}

If we consider solutions to equations without a gradient potential, a slightly modified proof leads to a better order of vanishing in the setting where $t \in \pb{1,\infty}$.
Although vanishing order estimates are not explicitly stated in \cite{KT16}, such results follow from the quantitative uniqueness theorems presented in that paper. 
The following theorem improves on the estimates implied from \cite{KT16} in two ways: we reduce the vanishing order and we extend the range of $t$ from $t > 2$ to all admissible exponents, $t > 1$.

\begin{theorem}
Assume that $\norm{V}_{L^t\pr{B_{R_0}}} \le M$ for some $t \in \pb{1, \iny}$.
Let $u$ be a solution to $\LP u + V u = 0$ in $B_{10}$.
Assume that $u$ is bounded and normalized in the sense of \eqref{bound} and \eqref{normal}. Then for any sufficiently small $\eps > 0$, the vanishing order of $u$ in $B_{1}$ is less than $C\pr{1 + C_2M^\mu}$.
That is, for any $x_0\in B_1$ and every $r$ sufficiently small,
\begin{align*}
 \|u\|_{L^\iny(B_{r}(x_0))}
 &\ge c r^{C\pr{1 + C_2M^\mu}},
\end{align*}
where
$\disp \mu = \left\{\begin{array}{ll}
\frac{2t}{3t-2} & 2 < t \le \iny \medskip \\
\frac{t}{2t-2-\eps\pr{2t-1-2\eps}} & 1 < t \le 2
\end{array}\right.$,
$c = c\pr{t, M, \hat C, \eps}$,
and $C_2 = C_2\pr{t,\eps}$. \label{thhh}
\end{theorem}

\begin{remark}
If we tried to prove Theorem \ref{thhh} using the same approach that gave Corollary \ref{thhW}, the resulting theorem would be weaker.
That is, if $W \equiv 0$, then Theorem \ref{thh} with $s = \iny$ and $K = 0$ implies a version of Theorem \ref{thhh} with $\mu$ replaced by
$\disp \tilde \mu = \left\{\begin{array}{ll}
\frac{2}{3} & t = \iny \medskip \\
\frac{2t}{3t-4-\eps(2t-4)}  & 2 < t < \iny \medskip \\
\frac{t}{t-1+\eps\pr{t-2\eps t}} & 1 < t \le 2
\end{array}\right.$.
For all $t \in \pr{1, \iny}$, we see that $\tilde \mu > \mu$, and therefore, such an approach gives a worse result than the one presented in Theorem \ref{thhh}.

This observation implies that the gradient potential $W$ plays some role in how the vanishing order depends on the norm of $V$. 
Notice that the different cases in Theorem \ref{thh} are determined by the relationship between $t$ and $s$, and that, in some cases, $\kappa$ and $\mu$ depends on both $t$ and  $s$.
Where these relationships comes from can be seen in the proofs below, Lemma \ref{CarlpqVW} for example.
\end{remark}

Finally, as a consequence of Theorem \ref{thhh}, we derive a quantitative unique continuation at infinity theorem.

\begin{theorem}
Assume that $\norm{V}_{L^t\pr{\R^2}} \le A_0$ for some $t \in \pb{1, \iny}$.
Let $u$ be a solution to $\LP u + V u = 0$ in $\R^2$.
Assume that $\norm{u}_{L^\iny\pr{\R^2}} \le C_0$ and $\abs{u\pr{0}} \ge 1$.
Then for any sufficiently small $\eps > 0$ and $R$ sufficiently large,
\begin{equation*}
\inf_{\abs{x_0} = R} \norm{u}_{L^\iny\pr{B_1\pr{x_0}}} \ge \exp\pr{-C R^\Pi \log R},
\end{equation*}
where
$\disp \Pi =  \left\{\begin{array}{ll}
\frac{4t-4}{3t-2} & t > 2 \medskip \\
\frac{2t-2}{2t-2-\eps\pr{2t-1-2\eps}} & 1 < t \le 2,
\end{array}\right.$,
and $C = C\pr{t, C_0, A_0, \eps}$.
\label{UCV}
\end{theorem}

Let's review some literature about unique continuation results in $\R^2$. 
The (weak) unique continuation property implies that a solution is trivial if the solution vanishes in an open subset in the domain. 
And the strong unique continuation property implies that a solution is trivial if the solution vanishes to infinite order at some point in the domain. 
The results of Schechter and Simon \cite{SS80} as well as Amrein, Berthier and Georgescu \cite{ABG81} show that solutions to \eqref{schro} in $\Om \su \R^2$ satisfy the weak unique continuation property whenever $V \in L^t_{loc}\pr{\Om}$, $t > 1$. 
Jerison and Kenig \cite{JK85} established the strong unique continuation property for (\ref{schro}) if $V \in L^{n/2}_{loc}\pr{\Om}$ for $n \ge 3$. 
In the setting where $n = 2$, they showed in \cite{JK85} that strong unique continuation holds for $V\in L^t_{loc}\pr{\Om}$ with $t > 1$.
On the other hand, the counterexample of Kenig and Nadirashvili \cite{KN00} implies that weak unique continuation can fail for $V \in L^1$. 
For drift equations of the form $\LP u + W \cdot \gr u = 0$ in $\Om \su \R^2$, the
counterexamples due to Mandache \cite{Man02} and Koch and Tataru \cite{KT02} show that weak unique continuation can fail for $W \in L^{2^-}$ and $L^2_{weak}$, respectively. 
It has been believed for some time (see the comments in \cite{KT02}, for example) that strong
unique continuation holds for $W \in L^{2^+}$. 
In the setting where $V$ and $W$ are assumed to be real-valued, with $s > 2$ and $t > 1$,
Allesandrini \cite{Ale12} used the correspondence between elliptic equations in the plane and first-order Beltrami equations to prove the strong unique continuation property. 
In fact, his results are more general since the leading operator may be variable with non-smooth, non-symmetric coefficients. 
Since we are able to characterize the vanishing order of solutions to \eqref{goal} with $V \in L^{1^+}$ and $W \in L^{2^+}$, in some sense, our results provide a complete description of quantitative uniqueness for elliptic equations in $\mathbb R^2$.

As in \cite{DZ17}, our main tool is an $L^p - L^q$ Carleman estimate.
To prove such an inequality, we replace Sogge's eigenfunction estimates \cite{Sog86} on $S^{n-1}$ for $n \ge 3$ with eigenfunction estimates that we derive explicitly from Parseval's inequality and from the fact that all eigenfunctions on $S^1$ are bounded.
The resulting Carleman estimates hold for $1< p \le 2 < q < \infty$.
Moreover, compared with corresponding results in \cite{DZ17},  the exponent of the parameter $\tau$ in the $n = 2$ Carleman estimate is always positive for all $p$ and $q$ within the appropriate ranges.
See Theorem \ref{Carlpq} for the details.
As a consequence, we can treat all $s > 2$ and all $t > 1$, thereby leaving no gaps between our results and we would expect from the literature.

The outline of the paper is as follows.
In section \ref{CarlEst}, we state and prove the main $L^q - L^q$ Carleman estimates.
Section \ref{Thm1Proof} is devoted to the proof of Theorem \ref{thh}.
The proof of Theorem \ref{thhh} is very similar to that of Theorem \ref{thh}, but the main differences are described in section \ref{Thm3Proof}.
In section \ref{QuantUC}, we show how the quantitative unique continuation at infinity theorems follow from the vanishing order estimates through a scaling argument.
Finally, we present the proof of an important lemma of $L^p\to L^q$ Carleman estimates and the quantitative Caccioppoli inequality in the appendix.
The letters $c$, $C$, $C_0$, $C_1$ and $C_2$ are independent of $u$ and may vary from line to line.

\section{Carleman estimates }
\label{CarlEst}

In this section, we state and prove the crucial tools, the $L^p-L^{ q}$ type Carleman estimates.
Let $r=|x-x_0|$.
Define the weight function to be
$$\phi(r)=\log r+\log(\log r)^2.$$
We use the notation $\|u\|_{L^p(r^{-2} dx)}$ to denote the $L^p$ norm with weight $r^{-2}$, i.e. $\|u\|_{L^p(r^{-2} dx)} = \disp \pr{\int |u\pr{x}|^p r^{-2}\, dx}^\frac{1}{p}$.
Our $L^p - L^q$ Carleman estimate for the Laplacian is as follows.

\begin{theorem}
Let $1< p \le 2 < q < \infty$.
For any $\eps \in \pr{0,1}$, there exists a constant $C$ and sufficiently small $R_0$ such that for any $u\in C^{\infty}_{0}\pr{B_{R_0}(x_0)\backslash\set{x_0} }$ and $\tau>1$,
one has
\begin{align}
&\tau^{1+ \be_1} \|(\log r)^{-1} e^{-\tau \phi(r)}u\|_{L^2(r^{-2}dx)}
+\tau^{\be_0} \|(\log r)^{-1} e^{-\tau \phi(r)}u\|_{L^q(r^{-2\pr{1 - \eps}}dx)}
\nonumber
\\ &+ \tau^{\be_1} \|(\log r )^{-1} e^{-\tau \phi(r)}r \nabla
u\|_{L^2(r^{-2}dx)}
\leq  C \|(\log r ) e^{-\tau \phi(r)} r^2 \LP u\|_{L^p(r^{-2}
dx)} ,
\label{mainCar}
\end{align}
where $\be_0 = \frac{2}{q}\pr{1 - \eps}+1-\frac{1}{p}$ and $\be_1= 1- \frac{1}{p}$.
Furthermore, $C = C\pr{p, q, \eps}$.
 \label{Carlpq}
\end{theorem}

To prove our Carleman estimate, we first establish some intermediate Carleman estimates for first-order operators.
Towards this goal, we introduce polar coordinates in $\mathbb R^2\backslash \{0\}$ by setting $x=r\omega$, with $r=|x|$ and $\omega=(\omega_1, \omega_2)\in S^{1}$.
Define a new coordinate $t=\log r$.
Then
$\disp \frac{\partial }{\partial x_j}=e^{-t}(\omega_j\partial_t+  \Omega_j)$, for $j = 1,2$,
where each $\Omega_j$ is a vector field in $S^{1}$.
It is well known that the vector fields $\Omega_j$ satisfy
$$\disp \sum^{2}_{j=1}\omega_j\Omega_j=0 \quad{and} \disp \quad
\sum^{2}_{j=1}\Omega_j\omega_j=1.$$
In the new coordinate system, the Laplace operator takes the form
\begin{equation}
e^{2t} \LP=\partial^2_t u+\LP_{\omega}, \label{laplace}
\end{equation}
where $\disp \LP_\omega=\sum_{j=1}^2 \Omega^2_j$ is the
Laplace-Beltrami operator on $S^{1}$.

The eigenvalues for $-\LP_\omega$ are $k^2$, $k\in \Z_{\ge 0}$.
The corresponding eigenspace is $E_k$, the space of spherical harmonics of degree $k$.
It follows that
$$\| \LP_\omega v\|^2_{L^2(dtd\omega)}=\sum_{k\geq 0} k^4\| v_k\|^2_{L^2(dtd\omega)}$$
and
\begin{equation}
\sum_{j=1}^n\| \Omega_j v\|^2_{L^2(dtd\omega)} =\sum_{k\geq 0}
k^2\|v_k\|^2_{L^2(dtd\omega)}, \label{lll}
\end{equation}
where $v_k$ denotes the projection of $v$ onto $E_k$ and
$\|\cdot\|_{L^2(dtd\omega)}$ denotes the $L^2$ norm on $(-\infty,
\infty)\times S^{1}$. Note that the projection operator, $P_k$, acts
only on the angular variables. In particular, $P_k v\pr{t, \om} =
P_k v\pr{t, \cdot} \pr{\om}$.
 Let
$$\Lambda=\sqrt{-\LP_\omega}.$$
The operator $\Lambda$ is a first-order elliptic pseudodifferential operator on $L^2(S^{1})$.
The eigenvalues for the operator $\Lambda$ are $k$, with corresponding eigenspace $E_k$.
That is, for any $v\in C^\infty_0(S^{1})$,
\begin{equation}
\Lambda  v= \sum_{k\geq 0}k P_k v.
\label{ord}
\end{equation}
Set
\begin{equation} L^\pm=\partial_t\pm \Lambda.
\label{use}
\end{equation} From the equation (\ref{laplace}), it follows that
\begin{equation*}
e^{2t} \LP=L^+L^-=L^-L^+.
\end{equation*}

With $r=e^t$, we define the weight function in terms of $t$,
$$\varphi(t)=\phi(e^t)=t+\log t^2.$$
We only consider the solutions in balls with small radius $r$. In
term of $t$,  we study the case when $t$ is sufficiently close to
$-\infty$.

We first state an $L^2- L^2$ Carleman inequality for the operator $L^+$.
For the proof of this result (which still holds when $n = 2$), we refer the reader to our companion paper, \cite{DZ17}.

\begin{lemma}
If $\abs{t_0}$ is sufficiently large, then for any $v \in C^\iny_0\pr{\pr{-\iny, -\abs{t_0}} \times S^{1}}$, we have that
\begin{eqnarray}
\tau \norm{t^{-1} e^{-\tau \varphi(t)}v}_{L^2(dtd\omega )}
&+&\norm{t^{-1}e^{-\tau \varphi(t)} \partial_t v}_{L^2(dtd\omega )}
+\sum_{j=1}^2 \norm{t^{-1}e^{-\tau \varphi(t)} \Omega_j v }_{L^2(dtd\omega )}   \nonumber \medskip\\
&\leq& C\norm{t^{-1} e^{-\tau \varphi(t)} L^+v}_{L^2(dtd\omega )}.
\label{cond}
 \end{eqnarray}
\label{Car22}
\end{lemma}

To prove the $L^p -L^2$ Carleman estimates, we use estimates for the projections of functions into the spaces of spherical harmonics.
The eigenfunctions on $S^1$ are those elements of $L^2\pr{S^1}$ that satisfy the system
\begin{equation}
\left\{ \begin{array}{lll} -\LP_\te e_\lambda &=&\lambda
e_\lambda, \quad 0\leq
\theta\leq 2\pi \nonumber \\
e_\lambda(0)&=&e_\lambda(2\pi), \end{array} \right.
\end{equation}
where we now use $\LP_\te$ to denote the Laplace-Beltrami operator on $S^1$. The orthonormal eigenfunctions are given by
$\frac{1}{\sqrt{\pi}}\cos k\theta$ and $\frac{1}{\sqrt{\pi}}\sin
k\theta$ for $k\in \mathbb Z_{\geq 0}$ with eigenvalue $\la_k =
k^2$. We can describe them
 by $e_k\pr{\te} =
\frac{1}{\sqrt{2\pi}} e^{i k \te}$ with $k \in \Z$. The following
result is the $n = 2$ analog of Lemma 3 in \cite{DZ17}. Since
Sogge's estimates on $S^{n-1}$ from \cite{Sog86} hold in $n \ge 3$,
our proof instead relies on Parseval's inequalities and the fact
that every eigenfunction is bounded.

\begin{lemma}
Let $M, N \in \N$ and let $\set{c_k}$ be a sequence of numbers such that $\abs{c_k} \le 1$ for all $k$.
For any $v \in L^2\pr{S^{1}}$ and every $p \in \brac{1, 2}$, we have that
\begin{align}
\|\sum^M_{k=N} c_k P_k v\|_{L^2(S^{1})}
&\leq C \pr{ \sum^M_{k=N} |c_k|^2}^{\frac{1}{p}-\frac{1}{2}}\| v\|_{L^p(S^{1})}.
\label{haha}
\end{align}
\label{upDown}
\end{lemma}

\begin{proof}
Recall that $P_k v = v_k$ is the projection of $v$ on the space of spherical harmonics of degree $k$, $E_k$.
Thus, $P_0 v = \innp{v, e_0} e_0$ and $P_k v = \innp{v, e_k} e_k + \innp{v, e_{-k}} e_{- k}$ for every $k \in \Z_+$, where we use the notation $\innp{\cdot, \cdot}$ to denote a pairing of elements in dual spaces.
By Parseval's identity,
$$ \sum^\infty_{k=-\infty}|\innp{e_k, \ v}|^2  = \|v\|_{L^2(S^1)}^2.  $$
Since $\|e_k\|_{L^\infty (S^1)}= \frac{1}{\sqrt{2\pi}}$ for all $k \in \Z$, then for $k \ne 0$,
\begin{align*}
\|P_k v\|_{L^\infty (S^1)}^2
&= \| \innp{v, e_k} e_k + \innp{v, e_{-k}} e_{-k}\|_{L^\infty (S^1)}^2
\le 2 \| \innp{v, e_k} e_k\|_{L^\infty (S^1)}^2
+ 2 \|\innp{v, e_{-k}} e_{-k}\|_{L^\infty (S^1)}^2 \\
&= 2 |\innp{e_k, \ v}|^2 \|e_k\|_{L^\infty (S^1)}^2
+ 2 |\innp{e_{-k}, \ v}|^2 \|e_{-k}\|_{L^\infty (S^1)}^2
= \frac{1}{\pi} \pr{|\innp{e_k, \ v}|^2 + |\innp{e_{-k}, \ v}|^2} \\
&\le \frac{1}{\pi} \sum^\infty_{j=-\infty}|\innp{e_j, \ v}|^2 =
\frac{1}{\pi} \|v\|_{L^2(S^1)}^2.
\end{align*}
Similarly, $\|P_0 v\|_{L^\infty (S^1)}^2  \le \frac 1 {2\pi} \|v\|_{L^2(S^1)}^2$.
Thus, for every $k \in \Z_{\ge 0}$,
\begin{equation}
\|P_k v\|_{L^\infty (S^1)}\leq c \|v\|_{L^2(S^1)}. \label{dim2}
\end{equation}
It follows that for any $u\in L^2(S^1)$
$$ \innp{u, P_k v}  = \innp{P_k u, v} \le \| P_k u\|_{L^{\iny}(S^{1})}\|v\|_{L^{1}(S^{1})} \le c \| u\|_{L^{2}(S^{1})}\|v\|_{L^{1}(S^{1})}. $$
By duality, we conclude that
\begin{align}
\| P_k v\|_{L^2\pr{S^1}} \le c \| v\|_{L^{1}(S^{1})}.
\label{down}
\end{align}
Since each $e_k$ is normalized in $L^2\pr{S^1}$, then $\sqrt{2\pi}
\norm{e_k}_{L^\iny\pr{S^1}} =1 = \norm{e_k}_{L^2\pr{S^1}}$ and
therefore, for $k \ne 0$,
\begin{align*}
\| P_k v\|_{L^\infty\pr{S^1}}
&= \| \innp{v, e_k} e_k + \innp{v, e_{-k}} e_{-k}\|_{L^\infty\pr{S^1}}
\le \abs{\innp{v, e_k} } \| e_k\|_{L^\infty\pr{S^1}} + \abs{\innp{v, e_{-k}}} \| e_{-k}\|_{L^\infty\pr{S^1}} \\
&= \frac 1 {\sqrt{2 \pi}} \pr{ \abs{\innp{v, e_k} } \| e_k\|_{L^2\pr{S^1}} + \abs{\innp{v, e_{-k}}} \| e_{-k}\|_{L^2\pr{S^1}}} \\
&= \frac 1 {\sqrt{2 \pi}} \pr{ \|\innp{v, e_k} e_k+ \innp{v, e_{-k}} e_{-k}\|_{L^2\pr{S^1}}}
= \frac 1 {\sqrt{2 \pi}}  \| P_k v \|_{L^2\pr{S^1}},
\end{align*}
where we have used orthogonality.
Similarly, $\| P_0 v\|_{L^\infty\pr{S^1}} = \frac 1 {\sqrt{2 \pi}}  \| P_0 v \|_{L^2\pr{S^1}}$.
Thus, $\| P_k v\|_{L^\infty\pr{S^1}} \le \frac 1 {\sqrt{2\pi}} \| P_k v\|_{L^2\pr{S^1}}$ for every $k \in \Z_{\ge 0}$.
Combining this observation with \eqref{down} shows that
\begin{align}
\| P_k v\|_{L^\infty\pr{S^1}} \le C \| v\|_{L^{1}(S^{1})}.
\label{downs}
\end{align}
Finally, it follows from the normalization condition in combination with Parseval's identity that
\begin{equation}
 \| P_k v\|_{L^{2}(S^{1})} \leq \|v\|_{L^{2}(S^{1})}.
 \label{stay}
\end{equation}
Interpolating \eqref{down} and \eqref{stay} gives that
\begin{align}
\| P_k v\|_{L^{2}(S^{1})} &\leq C(p) \|v\|_{L^{p}(S^{1})}
\label{indu}
\end{align}
for all $1\leq p\leq 2$.

Now we consider a more general setting.
Let $\{c_k\}$ be a sequence of numbers with $|c_k|\leq 1$.
For all $M\leq N$, it follows from orthogonality and H\"older's inequality that
\begin{align*}
\|\sum^M_{k=N} c_k P_k v\|^2_{L^2(S^{1})}
= \sum^M_{k=N} \| c_k P_k v\|^2_{L^2(S^{1})}
= \sum^M_{k=N} \abs{c_k}^2 \innp{P_k v, v}
\leq \sum^M_{k=N} |c_k|^2 \|P_k v\|_{L^{\infty}(S^{1})}\|v\|_{L^{1}(S^{1})}.
\end{align*}
Combining this inequality with (\ref{downs}) shows that
\begin{equation*}
\|\sum^M_{k=N} c_k P_k v\|_{L^2(S^{1})} \leq C \pr{\sum^M_{k=N}
|c_k|^2}^{\frac{1}{2}} \|v\|_{L^{1}(S^{1})}.
\end{equation*}
Clearly, as long as $\abs{c_k} \le 1$, then
\begin{equation}
\|\sum^M_{k=N} c_k P_k v\|_{L^2(S^{1})} \leq \| v\|_{L^2(S^{1})}.
\label{seqSame}
\end{equation}
As before, we interpolate the last two inequalities to reach \eqref{haha}.
\end{proof}

Now we state an $L^p-L^2$ type Carleman estimate for the operator $L^-$.

\begin{lemma}
For every $v \in C^\iny_c\pr{(-\infty, \ t_0)\times S^{1}}$ and $1<
p \le 2$,
\begin{equation}
\|t^{-{1}} e^{-\tau \varphi(t)} v\|_{L^2(dtd\omega)} \leq C
\tau^\beta \|t e^{-\tau \varphi(t)} L^- v\|_{L^p(dtd\omega)},
\label{key-}
\end{equation}
where $\be =  \frac{1-p}{p}$.
\label{CarLpp}
\end{lemma}

When $p = 2$, the proof of Lemma \ref{CarLpp} follows from the proof of Lemma 2 in \cite{DZ17} (which still holds when $n = 2$) combined with the fact that $\abs{t} \ge \abs{t_0} \ge C$. 
For $p \in \pr{1, 2}$, the proof of Lemma \ref{CarLpp} is similar to the proof of Lemma 4 in \cite{DZ17}. 
The major difference is that we need to replace the eigenfunction estimates from Lemma 3 in \cite{DZ17} by Lemma \ref{upDown} in the current paper. 
For the readers' convenience and the complete the presentation, we include the proof of Lemma \ref{CarLpp} in the Appendix.

Lemma \ref{CarLpp} plays the crucial role in the $L^p\to L^q$ Carleman estimates.
Let's compare the parameter $\beta$ from Lemma 4 in \cite{DZ17} with the one in Lemma \ref{CarLpp} above.
In Lemma \ref{CarLpp} above, we see that $\beta<0$ for any $p>1$.
However, in Lemma 4 from \cite{DZ17}, $\beta<0$ if $p>\frac{6n-4}{3n+2}$, so $\beta$ changes signs over the full range of $p$ values.
The fact that $\be$ doesn't change signs here implies that we can deal with all the admissible $s$
and $t$ in Theorem \ref{thh}.
The proof of Lemma 4 in \cite{DZ17} was partially motivated by the work in \cite{Reg99}, \cite{Jer86} and \cite{BKRS88}.

We now have all of the intermediate results required to prove the general $L^p- L^q$ Carleman  given in Theorem \ref{Carlpq}.
We combine Lemmas \ref{Car22} and \ref{CarLpp}, apply a Sobolev inequality, then interpolate to reach the conclusion of Theorem \ref{Carlpq}.

\begin{proof}[Proof of Theorem \ref{Carlpq}]
Let $u\in C^{\infty}_{0}\pr{B_{R_0}(x_0)\backslash\set{x_0} }$.
After to an application of Lemma \ref{Car22} to $v$, we apply Lemma \ref{CarLpp} to $L^+v$, and see that
\begin{align*}
\tau \norm{t^{-1} e^{-\tau \varphi(t)}v}_{L^2(dtd\omega )}
&+\norm{t^{-1}e^{-\tau \varphi(t)} \partial_t v}_{L^2(dtd\omega )}
+\sum_{j=1}^2 \norm{t^{-1}e^{-\tau \varphi(t)} \Omega_j v }_{L^2(dtd\omega )}   \\
&\leq C\norm{t^{-1} e^{-\tau \varphi(t)} L^+v}_{L^2(dtd\omega )}
\le C\tau^{\be}\norm{t e^{-\tau \varphi(t)} L^-
L^+v}_{L^p(dtd\omega )},
\end{align*}
where $\be = \frac{1-p}{p}$.
Recalling the definitions of $t$, $\vp$, and $L^\pm$, this gives
\begin{align}
&\tau \|(\log r)^{-1} e^{-\tau \phi(r)}u\|_{L^2(r^{-2}dx)}
+ \|(\log r )^{-1} e^{-\tau \phi(r)}r \nabla u\|_{L^2(r^{-2}dx)}
\nonumber \medskip\\
&\leq  C \tau^{\beta}  \|(\log r ) e^{-\tau \phi(r)} r^2 \LP u\|_{L^p(r^{-2}dx)}.
\label{dodo}
\end{align}
Direct computations shows that $\phi'(r)=\frac{1}{r}+\frac{2}{r\log
r} \le \frac 1 r$ since $r \le R_0 \le 1$. By the Sobolev imbedding
$W^{1, 2}\hookrightarrow L^{q'}$ with any $2<q < q' <\infty$ in
$n=2$, we have
\begin{align}
 \|(\log r)^{-1} e^{-\tau \phi(r)} u\|_{L^{q'}}
&\le c_{q'} \|\nabla [(\log r)^{-1} e^{-\tau \phi(r)}  u ]\|_{L^2} \nonumber \\
&\le C\tau \| (\log r)^{-1} e^{-\tau \phi(r)} r^{-1}  u\|_{L^2}
+C \| (\log r)^{-1} e^{-\tau \phi(r)}  \nabla u \|_{L^2}  \nonumber \\
&+C \| (\log r)^{-2} e^{-\tau \phi(r)} r^{-1} u \|_{L^2} \nonumber \\
&\le C\tau \| (\log r)^{-1} e^{-\tau \phi(r)} u \|_{L^2(r^{-2}dx)}
+C \| (\log r)^{-1} e^{-\tau \phi(r)} r \nabla u \|_{L^2(r^{-2}dx)}  \nonumber \\
&+C \| (\log r)^{-2} e^{-\tau \phi(r)}  u \|_{L^2(r^{-2}dx)} \nonumber \\
&\leq  C \tau^\be \|(\log r ) e^{-\tau \phi(r)} r^2 \LP
u\|_{L^p(r^{-2} dx)}, \label{L2*Est}
\end{align}
where the last inequality follows from \eqref{dodo}. Obviously, the
inequality \eqref{dodo} indicates that
\begin{equation}
\|(\log r)^{-1} e^{-\tau \phi(r)}u\|_{L^2\pr{r^{-2}dx}} \leq  C \tau^{\be -1}
\|(\log r ) e^{-\tau \phi(r)} r^2 \LP u\|_{L^p(r^{-2} dx)}.
\label{L2Est}
\end{equation}
To get a range of $L^q$-norms on the left, we interpolate the last two inequalities.
Choose $\la \in \pr{0,1}$ so that $q = 2 \la + \pr{1-\la} q'$.
The application of H\"older's inequality yields that
\begin{align*}
\|(\log r)^{-1} e^{-\tau \phi(r)}u\|_{L^q\pr{r^{-2\la}dx}}
&\le \|(\log r)^{-1} e^{-\tau \phi(r)}u\|_{L^2\pr{r^{-2}dx}}^{\frac{2\la}{q}} \|(\log r)^{-1} e^{-\tau \phi(r)}u\|_{L^{q'}}^{\frac{q'\pr{1-\la}}{q}}.
\end{align*}
Since $\la  = \frac{q'-q}{q' -2}$, if we set $\te = \frac{2\la}{q} = \frac{2(q'-q)}{q(q'-2)}$, then $1 - \te = 1-\frac{2(q'-q)}{q(q'-2)}= \frac{q'(q-2)}{q(q'-2)}$.
Thus, $\theta$ falls in the interval $[0,\ 1]$.
From \eqref{dodo} and \eqref{L2Est}, we obtain that
\begin{align*}
& \|(\log r)^{-1} e^{-\tau \phi(r)}u\|_{L^q\pr{r^{-2\la}dx}} \\
&\leq  \|(\log r)^{-1} e^{-\tau \phi(r)}u\|_{L^2\pr{r^{-2}dx}}^{\theta}
\|(\log
r)^{-1} e^{-\tau\phi(r)}u\|_{L^{q'}}^{1-\theta} \\
&\le \brac{C\tau^{\beta-1} \|(\log r ) e^{-\tau \phi(r)} r^2
\LP u\|_{L^p(r^{-2} dx)}}^\te
\brac{C \tau^\beta  \|(\log r ) e^{-\tau \phi(r)} r^2 \LP u\|_{L^p(r^{-2} dx)}}^{1 - \te} \\
&= C \tau^{\beta-\te} \|(\log r ) e^{-\tau \phi(r)} r^2 \LP
u\|_{L^p(r^{-2} dx)}.
\end{align*}
We have $\be - \te = \frac{1}{p}-1-\frac{2}{q}\pr{1-\eps}$, where $\eps=\frac{q-2}{q'-2}$.
Since $q' \in \pr{q,\infty}$ and $q > 2$, then $\eps \in \pr{0,1}$. Moreover, since we can choose $q'$ to be arbitrarily large, then $\eps$ may be made arbitrarily close to zero.
Recalling the definition of $\la$, we conclude that for any $2 < q<\infty$ and any $\eps \in \pr{0,1}$,
\begin{equation}
\tau^{\frac{2}{q}\pr{1 - \eps}+1-\frac{1}{p}}\|(\log r)^{-1} e^{-\tau\phi(r)}u\|_{L^q\pr{r^{-2\pr{1-\eps}}dx}}
\leq C \|(\log r ) e^{-\tau \phi(r)} r^2 \LP u\|_{L^p(r^{-2} dx)}.
\label{inter}
\end{equation}
Combining \eqref{inter} with \eqref{dodo}, we arrive at the proof of Theorem \ref{Carlpq}.
\end{proof}

\section{The proof of Theorem \ref{thh}}
\label{Thm1Proof}

The first step in the proof of Theorem \ref{thh} is to establish a Carleman estimate for the operator $\LP + W \cdot \gr + V$.
We use the triangle inequality and H\"older's inequality along with the crucial Carleman estimates in Theorem \ref{Carlpq}.
Because of the correlation between the potentials $W(x)$ and $V(x)$, we need to work in cases depending on the relationships between $s$ and $t$.

\begin{lemma}
Assume that for some $s \in \pb{2, \ \iny}$ and $t \in \pb{ 1, \ \infty}$, $\norm{W}_{L^s\pr{B_{R_0}}} \le K$ and $\norm{V}_{L^t\pr{B_{R_0}}} \le M$.
Then for every sufficiently small $\eps > 0$, there exist constants $C_0$, $C_1$, $C_2$, and sufficiently small $R_0 < 1$  such that for
any $u\in C^{\infty}_{0}(B_{R_0}(x_0)\setminus \set{x_0})$ and
$$\tau \ge 1+ C_1 K^{\kappa} + C_2 M^{\mu},$$
one has
\begin{align}
\tau^{2 - \frac 1 p} \|(\log r)^{-1} e^{-\tau \phi(r)}u\|_{L^2(r^{-2}dx)}
&\leq  C_0 \|(\log r ) e^{-\tau \phi(r)} r^2\pr{ \LP u + W
\cdot \gr u + V u}\|_{L^p(r^{-2} dx)} , \label{main1}
\end{align}
where
$\disp \kappa = \left\{\begin{array}{ll}
\frac{2s}{s-2} & t > \frac{2s}{s+2} \medskip \\
\frac{t}{t-1-\eps t} & 1 < t \le \frac{2s}{s+2}
\end{array}\right.$,
$\disp \mu = \left\{\begin{array}{ll}
\frac{2s}{3s-2} & t \ge s \medskip \\
\frac{2st}{3st+2t-4s-\eps(2st+4t-4s)}  & \frac{2s}{s+2} < t < s \medskip \\
\frac{t}{t-1+\eps\pr{t-2\eps t}} & 1 < t \le \frac{2s}{s+2}
\end{array}\right.$,
 and \newline
 $\disp p = \left\{\begin{array}{ll}
\frac{2s}{s+2} & t > \frac{2s}{s+2} \medskip \\
\frac{t}{1+\eps t} & 1 < t \le \frac{2s}{s+2}
\end{array}\right.$.
Moreover, $C_0 = 2C$, where $C\pr{p, q, \eps}$ is from Theorem \ref{Carlpq} with $\disp q =
\left\{\begin{array}{ll}
\frac{2st}{st+2t-2s} & \frac{2s}{s+2} < t < s \medskip \\
\frac 1 \eps & 1 < t \le \frac{2s}{s+2}
\end{array}\right.$,
$C_1 = C_1\pr{s, t, \eps}$,
and $C_2 = C_2\pr{s,t, \eps}$.
\label{CarlpqVW}
\end{lemma}

\begin{proof}
Assume that $\eps > 0$ is sufficiently small.
From our crucial estimate \eqref{mainCar} in Theorem \ref{Carlpq} and the triangle inequality, we get that
\begin{align}
&\tau^{2 - \frac 1 p} \|(\log r)^{-1} e^{-\tau \phi(r)}u\|_{L^2(r^{-2}dx)}
+\tau^{\frac 2 q \pr{1 - \eps} + 1 - \frac 1 p} \|(\log r)^{-1} e^{-\tau \phi(r)}u\|_{L^q(r^{-2\pr{1 - \eps}}dx)} \nonumber \\
&+\tau^{1 - \frac 1 p} \|(\log r )^{-1} e^{-\tau \phi(r)}r \nabla u\|_{L^2(r^{-2}dx)} \nonumber \\
&\le C\|(\log r) e^{-\tau \phi(r)} r^2 (\LP u)\|_{L^p(r^{-2} dx)} \nonumber \\
&\le C\|(\log r) e^{-\tau \phi(r)} r^2 (\LP u+ W\cdot \nabla u + V u)\|_{L^p(r^{-2} dx)} \nonumber \\
&+ C\|(\log r) e^{-\tau \phi(r)} r^2 W\cdot \nabla u \|_{L^p(r^{-2}
dx)} + C\|(\log r) e^{-\tau \phi(r)} r^2 V u \|_{L^p(r^{-2} dx)}.
\label{triIneq}
\end{align}
To reach the estimate (\ref{main1}) in the lemma, we will absorb the last two terms above into the lefthand side by appropriately choosing $p$ and $q$, and making $\tau$ sufficiently large.

If $p \in \pb{1, 2}$, by H\"older' inequality, it follows that
\begin{align}
& \|(\log r) e^{-\tau\phi(r)} r^2 W\cdot\nabla u\|_{L^p(r^{-2} dx)} \nonumber \\
&\le  \|W\|_{L^{\frac{2p}{2-p}}\pr{B_{R_0}}} \|\pr{\log r}^2 r^{2 - \frac 2 p}\|_{L^{\iny}\pr{B_{R_0}}}
\|(\log r)^{-1} e^{-\tau \phi(r)} r |\nabla u|\|_{L^2(r^{-2} dx)} \nonumber \\
&\le c \|W\|_{L^{\frac{2p}{2-p}}\pr{B_{R_0}}} \|(\log r)^{-1}
e^{-\tau \phi(r)} r |\nabla u|\|_{L^2(r^{-2} dx)},
\label{hod2}
\end{align}
where we have used the fact that $2 - \frac 2 p > 0$ and $R_0$ is small enough.
Similarly, an application of H\"older's inequality implies that
\begin{align}
\|(\log r) e^{-\tau \phi(r)} r^2 V u \|_{L^p(r^{-2} dx)}
&\le  \|V\|_{L^{\frac{2p}{2-p}}\pr{B_{R_0}}} \|\pr{\log r}^2 r^{3 - \frac 2 p}\|_{L^{\iny}\pr{B_{R_0}}}
\|(\log r)^{-1} e^{-\tau \phi(r)} u\|_{L^2(r^{-2} dx)} \nonumber \\
&\le c \|V\|_{L^{\frac{2p}{2-p}}\pr{B_{R_0}}} \|(\log r)^{-1}
e^{-\tau \phi(r)} u\|_{L^2(r^{-2} dx)}.
\label{hod3}
\end{align}
Note that $1 - \frac 1 p + \frac 1 q\pr{1 - \eps} > 0$.
Therefore, if $p \in \pb{1, 2}$ and $q > 2$, by H\"older's inequality again, we obtain that
\begin{align}
& \|(\log r) e^{-\tau \phi(r)} r^2 V u \|_{L^p (r^{-2}dx)} \nonumber \\
&\le  \|V\|_{L^{\frac{pq}{q-p}}\pr{B_{R_0}}} \|\pr{\log r}^2
r^{2\brac{1 - \frac 1 p + \frac 1 q\pr{1 - \eps}}}
\|_{L^{\iny}\pr{B_{R_0}}}
\|(\log r)^{-1} e^{-\tau \phi(r)} u r^{ -\frac 2 q\pr{1 - \eps}}\|_{L^q} \nonumber \\
&\le c \|V\|_{L^{\frac{pq}{q-p}}\pr{B_{R_0}}} \|(\log r)^{-1}
e^{-\tau \phi(r)} u\|_{L^q(r^{-2\pr{1 - \eps}}dx)}. \label{hod33}
\end{align}
Next we work in cases to achieve the conclusion (\ref{main1}) in the lemma and determine the appropriate power of $\tau$.
\\

\nid {\bf Case 1: $t \in \brac{s, \iny}$} \\
If $t \ge s$, then we choose $p = \frac{2s}{s+2}$.
Since $s > 2$, then $p\in (1, \ 2]$ is in the appropriate range of Theorem \ref{Carlpq}.
As $\frac{2p}{2-p} = s \le t$, by H\"older's inequality, $\|V\|_{L^{s}} \le c \|V\|_{L^{t}}$.
Substituting \eqref{hod2} and \eqref{hod3} into \eqref{triIneq} yields that
\begin{align*}
&\tau^{\frac 3 2 - \frac 1 s} \|(\log r)^{-1} e^{-\tau \phi(r)}u\|_{L^2(r^{-2}dx)}
+\tau^{\frac 1 2 - \frac 1 s} \|(\log r )^{-1} e^{-\tau \phi(r)}r \nabla u\|_{L^2(r^{-2}dx)} \nonumber \\
&\le C\|(\log r) e^{-\tau \phi(r)} r^2 (\LP u+ W\cdot \nabla u + V u)\|_{L^p(r^{-2} dx)} \nonumber \\
&+ c C K \|(\log r)^{-1} e^{-\tau \phi(r)} r |\nabla u|\|_{L^2(r^{-2} dx)}
+ c C M \|(\log r)^{-1} e^{-\tau \phi(r)} u\|_{L^2(r^{-2} dx)}.
\end{align*}
In order to  absorb the last two terms on the right into the lefthandside in the last inequality, we choose $\tau \ge 1 + \pr{c C K}^{\frac {2s} {s-2} }+ \pr{2c C M}^{\frac {2s} {3s-2}} $ to get the
conclusion (\ref{main1}).
\\

\nid {\bf Case 2: $t \in \pr{\frac{2s}{s+2}, s}$} \\
In this case, we choose $p = \frac{2s}{s+2}$.
We use (\ref{hod33}) to absorb the term involving the potential $V(x)$.
We need to choose $t=\frac{pq}{q-p}$.
Since $p = \frac{2s}{s+2}$, then $q=\frac{2st}{st+2t-2s}$.
We can check that $p$ falls in the appropriate range and that $q \in \pr{2, \infty}$ follows from the assumption on $t$.
Since $\frac{2p}{2-p} = s$ and $\frac{pq}{q-p} = t$, substituting \eqref{hod2} and \eqref{hod33} into
\eqref{triIneq}, we obtain that
\begin{align}
&\tau^{\frac 3 2 - \frac 1 s} \|(\log r)^{-1} e^{-\tau
\phi(r)}u\|_{L^2(r^{-2}dx)} +\tau^{\frac 3 2 + \frac 1 s - \frac 2 t
- \eps \pr{ 1 + \frac 2 s - \frac 2 t }}
\|(\log r)^{-1} e^{-\tau \phi(r)}u\|_{L^q(r^{-2\pr{1 - \eps}}dx)} \nonumber \\
&+\tau^{\frac 1 2 - \frac 1 s} \|(\log r )^{-1} e^{-\tau \phi(r)}r \nabla u\|_{L^2(r^{-2}dx)} \nonumber \\
&\le C\|(\log r) e^{-\tau \phi(r)} r^2 (\LP u+ W\cdot \nabla u + V u)\|_{L^p(r^{-2} dx)} \nonumber \\
&+ c C K \|(\log r)^{-1} e^{-\tau \phi(r)} r |\nabla
u|\|_{L^2(r^{-2} dx)} + c C M \|(\log r)^{-1} e^{-\tau \phi(r)}
u\|_{L^q(r^{-2\pr{1 - \eps}}dx)}.
\label{boundSub}
\end{align}
Since $s > 2$ implies that $3st + 2t - 4s>2st + 4t - 4s$, then $0<\eps<1<\frac{3st + 2t - 4s}{2st + 4t -4s}$ and we have that $\brac{\frac 3 2 + \frac 1 s - \frac 2 t - \eps\pr{ 1 + \frac 2 s - \frac 2 t }}^{-1}>0$.
Straightforward computations show that $\brac{\frac 3 2 + \frac 1 s - \frac 2 t - \eps\pr{ 1 + \frac 2 s - \frac 2 t }}^{-1} =\frac{2st}{3st+2t-4s-\eps ( 2st+4t-4s)}$.
Therefore, to absorb the last two terms into the lefthand side and get (\ref{main1}), we take $\tau \ge
1 + \pr{c C K}^{\frac {2s} {s-2}}+ \pr{c C M}^{\frac{2st}{3st+2t-4s-\eps ( 2st+4t-4s)}}$.
\\

\nid {\bf Case 3: $t \in \pb{1, \frac{2s}{s+2}}$} \\
As in the previous case, we'll substitute \eqref{hod2} and \eqref{hod33} into \eqref{triIneq}.
Therefore, we need to choose $p \in \pr{1,2}$ and $q \in \pr{2, \iny}$ so that $\frac{2p}{2-p} \le s$ and $\frac{pq}{q-p} = t$.
Let $q = \frac{1}{\eps}$ and $p =\frac{tq}{q+t}$.
Note that $p=\frac{tq}{q+t}=\frac{t}{1+\eps t} \in \pr{1, 2}$ if $0<\eps<\frac{t-1}{t}$.
And as long as $\eps < \frac 1 2$, then $q$ is in the appropriate range.
Since $1<t\leq \frac{2s}{s+2}$, then $s\geq \frac{2t}{2-t}$.
It follows that
$$\frac{2p}{2-p} = \frac{2tq}{2q+2t-tq}=\frac{2t}{2-t+2t\eps}<\frac{2t}{2-t} \le s.$$
An application of H\"older's inequality shows that $\|W\|_{L^{\frac{2p}{2-p}}} \le c \|W\|_{L^{s}}$.
We have
$$\frac{2}{q}\pr{1 - \eps}+1-\frac{1}{p}= 2 \eps \pr{1 - \eps}+1-\frac{1+\frac t q}{t}=\frac{t-1}{t}+\eps_0 >0,$$
where $\eps_0 = \eps\pr{1-2\eps}$.
Substituting \eqref{hod2} and
\eqref{hod33} into \eqref{triIneq} gives that
\begin{align}
&\tau^{\frac{2t-1}{t}-\eps} \|(\log r)^{-1} e^{-\tau
\phi(r)}u\|_{L^2(r^{-2}dx)}
+\tau^{\frac{t-1}{t}+\eps_0} \|(\log r)^{-1} e^{-\tau \phi(r)}u\|_{L^q(r^{-2\pr{1 - \eps}}dx)} \nonumber \\
&+\tau^{\frac{t-1}{t}-\eps} \|(\log r )^{-1} e^{-\tau \phi(r)}r \nabla u\|_{L^2(r^{-2}dx)} \nonumber \\
&\le C\|(\log r) e^{-\tau \phi(r)} r^2 (\LP u+ W\cdot \nabla u + V u)\|_{L^p(r^{-2} dx)} \nonumber \\
&+ c C K \|(\log r)^{-1} e^{-\tau \phi(r)} r |\nabla
u|\|_{L^2(r^{-2} dx)} + c C M \|(\log r)^{-1} e^{-\tau \phi(r)}
u\|_{L^q(r^{-2\pr{1 - \eps}}dx)}. \label{boundSub}
\end{align}
If we choose $\tau \ge 1 + \pr{ c C K}^{\frac{t}{t-1-\eps t}} + \pr{ c C M}^{\frac{t}{t-1+\eps_0t}}$,  we will arrive at the conclusion (\ref{main1}).
\end{proof}

With the aid of the Carleman estimates in Lemma \ref{CarlpqVW}, we present the $L^\iny$ three-ball inequality that will serve as an important tool in the proof of Theorem \ref{thh}.

\begin{lemma}
Let $0 < r_0< r_1< R_1 < R_0$, where $R_0 < 1$ is sufficiently small.
Assume that for some $s \in \pb{2, \iny}$, $t \in \pb{1, \iny}$, $\|W\|_{L^s\pr{B_{R_0}}} \le K$ and $\|V\|_{L^t\pr{B_{R_0}}} \le M$.
Let $u$ be a solution to \eqref{goal}.
Then, for any sufficiently small $\eps, \de > 0$,
\begin{align}
\|u\|_{L^\infty \pr{B_{3r_1/4}}}
&\le C F_\de\pr{r_1} |\log r_1| \brac{ (K+|\log r_0|)F_\de\pr{r_0} \|u\|_{L^\iny(B_{2r_0})}}^{k_0} \nonumber \\
&\times \brac{(K+|\log R_1|)  F_\de\pr{R_1}\|u\|_{L^\iny(B_{R_1})}}^{1 - k_0} \nonumber \\
&+C F_\de\pr{r_1} \pr{\frac{R_1 }{r_1}} \pr{1 +\frac{|\log r_0|}{K}} \nonumber \\
&\times \exp\brac{\pr{1 + C_1 K^\kappa + C_2 M^\mu} \pr{\phi\pr{\frac{R_1}{2}}-\phi(r_0)}} \|u\|_{L^\iny(B_{2r_0})},
\label{three}
\end{align}
where $\disp k_0 =
\frac{\phi(\frac{R_1}{2})-\phi(r_1)}{\phi(\frac{R_1}{2})-\phi(r_0)}$,
$F_\de\pr{r} = (1 + r K^{\frac{s}{s-2}+\de} + r M^{\frac{t}{2t-2}+\de})$, and $\kappa$, $\mu$, $C_1$, and $C_2$ are as given in Lemma \ref{CarlpqVW}, and $C = C\pr{s, t, \eps}$.
\label{threeBall}
\end{lemma}

\begin{proof}
Fix the sufficiently small constants $\eps, \de > 0$.
Let $r_0< r_1< R_1$.
Choose a smooth function $\eta\in C^\infty_{0}(B_{R_0})$ with $B_{2R_1}\subset B_{R_0}$.
The standard notation $\brac{a,b}$ is used to denote a closed annulus with inner radius $a$ and outer radius
$b$. Let
$$D_1=\brac{\frac{3}{2}r_0, \frac{1}{2}R_1 }, \quad  \quad
D_2= \brac{r_0, \frac{3}{2}r_0}, \quad \quad
D_3=\brac{\frac{1}{2}R_1, \frac{3 }{4}R_1}.$$
Let $\eta=1$ on $D_1$ and $\eta=0$ on $[0, \ r_0]\cup \brac{\frac{3}{4}R_1, \ R_1}$.
It is easy to see that $|\nabla \eta|\leq \frac{C}{r_0}$ and $|\nabla^2\eta|\leq \frac{C}{r_0^2}$ on $D_2$.
Similarly, we have $|\nabla \eta|\leq \frac{C}{R_1}$ and $|\nabla^2 \eta|\leq\frac{C}{R_1^2}$ on $D_3$.

Since $u$ is a solution to \eqref{goal} in $B_{R_0}$, as discussed in the introduction, $u \in L^\iny\pr{B_{R_1}} \cap W^{1,2}\pr{B_{R_1}} \cap W^{2,p}\pr{B_{R_1}}$.
By regularization, the estimate in Lemma \ref{CarlpqVW} holds for $\eta u$.
Taking into account that $u$ is a solution to equation \eqref{goal} and substituting $\eta u$ into the Carleman estimate in Lemma \ref{CarlpqVW}, we get that whenever
$$\tau \ge 1+ C_1 K^{\kappa} + C_2 M^{\mu},$$
\begin{align*}
\tau^{2 - \frac 1 p} \|(\log r)^{-1} e^{-\tau \phi(r)} \eta u\|_{L^2(r^{-2}dx)}
&\leq  C_0 \|(\log r ) e^{-\tau \phi(r)} r^2\pr{ \LP \pr{\eta u} + W \cdot \gr\pr{\eta u} + V \eta u}\|_{L^p(r^{-2} dx)} \\
&=  C_0 \|(\log r ) e^{-\tau \phi(r)} r^2\pr{ \LP \eta \, u + 2 \gr \eta \cdot \gr u + W \cdot \gr \eta \, u }\|_{L^p(r^{-2} dx)},
\end{align*}
where $\kappa$, $\mu$, and $p$ are as given in the Lemma \ref{CarlpqVW}.
Thus, the following holds
\begin{equation}
\tau^{2 - \frac 1 p} \|(\log r)^{-1} e^{-\tau \phi(r)} u\|_{L^2(D_1,
r^{-2}dx )}\leq \mathcal{B}, \label{jjj}
\end{equation}
where
$$ \mathcal{B}:= C_0 \|(\log r) e^{-\tau \phi(r)} r^2 (\LP \eta \, u+W\cdot \nabla \eta \, u + 2\nabla \eta \cdot \nabla u)\|_{L^p(D_2\cup D_3, r^{-2} dx)}.$$
To bound $\mathcal{B}$, an application of H\"older's inequality yields that
\begin{align*}
\|(\log r) e^{-\tau \phi(r)} r^2 \LP \eta \, u\|_{L^p(D_2\cup D_3, r^{-2} dx)}
&\le \| \pr{\log r} r^{2} \LP \eta\|_{L^\iny\pr{D_2}} \| r^{- \frac 2 p}\|_{L^{\frac{2p}{2-p}}\pr{D_2}}\| e^{-\tau \phi(r)} u\|_{L^2(D_2)} \\
&+\| \pr{\log r} r^{2} \LP \eta\|_{L^\iny\pr{D_2}} \| r^{- \frac 2
p}\|_{L^{\frac{2p}{2-p}}\pr{D_2}}\| e^{-\tau \phi(r)} u\|_{L^2(D_2)}
\end{align*}
and
\begin{align*}
\|(\log r) e^{-\tau \phi(r)} r^2 \gr \eta \cdot \gr u\|_{L^p(D_2\cup D_3, r^{-2} dx)}
&\le \| \pr{\log r} r^{2} \gr \eta\|_{L^\iny\pr{D_2}} \| r^{- \frac 2 p}\|_{L^{\frac{2p}{2-p}}\pr{D_2}}\| e^{-\tau \phi(r)} \gr u\|_{L^2(D_2)} \\
&+\| \pr{\log r} r^{2} \gr \eta\|_{L^\iny\pr{D_3}} \| r^{- \frac 2 p}\|_{L^{\frac{2p}{2-p}}\pr{D_3}}\| e^{-\tau \phi(r)} \gr u\|_{L^2(D_3)}.
\end{align*}
Since $\frac{2p}{2-p} \le s$, arguments that are similar to those that appear in \eqref{hod2} show that
\begin{align*}
&\|(\log r) e^{-\tau\phi(r)} r^2 W\cdot\nabla \eta \, u\|_{L^p(r^{-2} dx, D_2 \cup D_3)} \\
&\le c \|W\|_{L^{s}\pr{D_2}} \norm{\gr \eta}_{L^\iny\pr{D_2}}\| e^{-\tau \phi(r)} r u \|_{L^2(D_2, r^{-2} dx)} \\
&+ c \|W\|_{L^{s}\pr{D_3}} \norm{\gr \eta}_{L^\iny\pr{D_3}}\| e^{-\tau \phi(r)} r u \|_{L^2(D_3, r^{-2} dx)} \\
&\le cK r_0^{-1}\|e^{-\tau \phi(r)} u\|_{{L^2(D_2)}}
+ cK R_1^{-1}\|e^{-\tau \phi(r)} u \|_{L^2(D_3)},
\end{align*}
where we have used the bounds on $\abs{\gr \eta}$.
From the estimates of $\eta$ in $D_2$ and $D_3$, it follows that
\begin{eqnarray*}
\mathcal{B} &\leq & C |\log r_0| r_0^{-1}\pr{\|e^{-\tau \phi(r)}
u\|_{{L^2(D_2)}} + r_0 \|e^{-\tau \phi(r)} \nabla u\|_{{L^2(D_2)}}}
+CK r_0^{-1}\|e^{-\tau \phi(r)}u\|_{{L^2(D_2)}} \medskip \\
&+& C|\log R_1|R_1^{-1}\pr{\|e^{-\tau \phi(r)} u\|_{{L^2(D_3)}}+ R_1
\|e^{-\tau \phi(r)} \nabla u\|_{{L^2(D_3)}}}
+CK R_1^{-1}\|e^{-\tau \phi(r)} u\|_{{L^2(D_3)}}.
\end{eqnarray*}
Therefore,
\begin{eqnarray*}
\mathcal{B} &\leq & C\pr{K + |\log r_0|} r_0^{-1} e^{-\tau
\phi(r_0)} \pr{\| u\|_{{L^2(D_2)}}+ r_0
\|\nabla u\|_{{L^2(D_2)}}}  \medskip \\
&+& C\pr{K+ |\log R_1|} R_1^{-1}e^{-\tau \phi\pr{\frac{R_1}{2}}} \pr{\|
u\|_{{L^2(D_3)}}+ R_1 \| \nabla u\|_{{L^2(D_3)}}},
\end{eqnarray*}
where we have used the fact that $e^{-\tau\phi(r)}$ is a decreasing function with respect to $r$.
To estimate the gradient term, $\nabla u$, we use the the Caccioppoli inequality (see Lemma \ref{CaccLem} in the appendix) to get that
\begin{equation*}
\|\nabla u\|_{L^2(D_2)}\leq
\frac{C}{r_0}F_\de\pr{r_0} \| u\|_{L^2\pr{B_{2r_0}\backslash B_{{r_0}/{2}}}}
\end{equation*}
and
\begin{equation*}
\|\nabla u\|_{L^2(D_3)}\leq
\frac{C}{R_1} F_\de\pr{R_1} \| u\|_{L^2\pr{B_{R_1}\backslash B_{{R_1}/{4}}}},
\end{equation*}
where we adopt the notation $F_\de\pr{r} = 1+r K^{\frac{s}{s-2} + \de}+r M^{\frac{t}{2t-2}+\de}$. Therefore,
\begin{eqnarray*}
\mathcal{B}
&\leq& C\pr{K+|\log r_0|} r_0^{-1}e^{-\tau \phi(r_0)}F_\de\pr{r_0}\|u\|_{L^2(B_{2r_0})}  \nonumber \medskip \\
&+& C\pr{K+|\log R_1|} {R_1}^{-1}e^{-\tau \phi\pr{\frac{R_1}{2}}}F_\de\pr{R_1} \|u\|_{L^2(B_{R_1})}.
\end{eqnarray*}
Introduce a new set $D_4=\{r\in D_1, \ r\leq r_1\}$.
From \eqref{jjj} and that $\tau \ge 1$ and $2 - \frac 1 p > 0$, it follows that
\begin{align*}
\| u\|_{L^2 (D_4)}
&\le \tau^{2 - \frac 1 p}\| u\|_{L^2 (D_4)}
\le \tau^{2 - \frac 1 p} \|e^{\tau\phi(r)}(\log r) r\|_{L^\iny\pr{D_4}} \| (\log r)^{-1} e^{-\tau\phi(r)} u\|_{L^2 (D_4, r^{-2}dx)}  \\
&\le e^{\tau \phi(r_1)} |\log r_1| r_1 \mathcal{B},
\end{align*}
where we have considered that $e^{\tau\phi(r)}(\log r) r$ is
increasing on $D_1$ for $R_0$ sufficiently small. Adding $\|u\|_{L^2
\pr{B_{3r_0/2}}}$ to both sides of the last inequality and using the
upper bound on $\mathcal{B}$ implies that
\begin{align*}
\| u\|_{L^2 (B_{r_1})}
&\le C |\log r_1|  \pr{K+|\log r_0|} \pr{\frac{r_1}{r_0}} e^{\tau \brac{\phi(r_1)- \phi(r_0)}} F_\de\pr{r_0}  \|u\|_{L^2(B_{2r_0})} \\
&+ C |\log r_1| \pr{K+|\log R_1|} \pr{\frac{r_1}{R_1}} e^{\tau\brac{\phi(r_1) - \phi\pr{\frac{R_1}{2}}}} F_\de\pr{R_1} \|u\|_{L^2(B_{R_1})}.
\end{align*}
For the ease of the presentation,  we define
$$U_1 =\|u\|_{L^2(B_{2r_0})}, \quad \quad U_2=\|u\|_{L^2(B_{R_1})},$$
$$A_1 = C |\log r_1|  \pr{K+|\log r_0|} \pr{\frac{r_1}{r_0}} F_\de\pr{r_0},$$
and
$$A_2 = C |\log r_1| \pr{K+|\log R_1|} \pr{\frac{r_1}{R_1}} F_\de\pr{R_1}.$$
Then the previous inequality simplifies to
\begin{eqnarray}
\| u\|_{L^2 (B_{r_1})}  &\leq& A_1 \brac{\frac{\exp \pr{\phi(r_1)}}{\exp\pr{ \phi(r_0)}}}^\tau U_1 + A_2 \brac{\frac{\exp\pr{ \phi(r_1)}}{\exp\pr{ \phi\pr{\frac{R_1}{2}}}}}^\tau U_2.
\label{D4est}
\end{eqnarray}
Introduce another parameter $k_0$ as
$$\frac{1}{k_0}=\frac{\phi(\frac{R_1}{2})-\phi(r_0)}{\phi(\frac{R_1}{2})-\phi(r_1)}.$$
Recall that $\phi(r)=\log r+\log (\log r)^2$.
If $r_1$ and $R_1$ are fixed, and $r_0\ll r_1$, i.e. $r_0$ is sufficiently small, then $\frac{1}{k_0}\simeq \log \frac{1}{r_0}$.
Let
$$\tau_1 =\frac{k_0}{\phi\pr{\frac{R_1}{2}}-\phi(r_1)}\log\pr{\frac{A_2{U}_2}{A_1 {U}_1}}.$$
If $\tau_1 \ge 1 + C_1 K^\kappa + C_2 M^\mu$, then the calculations performed above are valid with $\tau = \tau_1$.
Substituting $\tau_1$ into \eqref{D4est} gives that
\begin{eqnarray}
\| u\|_{L^2 (B_{r_1})}  &\leq& 2\pr{A_1 U_1}^{k_0}\pr{A_2 U_2}^{1 - k_0}.
\label{mix1}
\end{eqnarray}
Instead, if $\tau_1 < 1 + C_1 K^\kappa + C_2 M^\mu$, then
\begin{align*}
U_2
< \frac{A_1}{A_2} \exp\brac{\pr{1 + C_1 K^\kappa + C_2 M^\mu} \pr{\phi\pr{\frac{R_1}{2}}-\phi(r_0)}} U_1.
\end{align*}
The last inequality indicates that
\begin{equation}
\|u\|_{L^2 (B_{r_1})}
\le C \pr{\frac{R_1}{r_0}} \pr{1 +\frac{|\log r_0|}{K}} e^{\pr{1 + C_1 K^\kappa + C_2 M^\mu} \pr{\phi\pr{\frac{R_1}{2}}-\phi(r_0)}} \|u\|_{L^2(B_{2r_0})}.
\label{mix2}
\end{equation}
The combination of \eqref{mix1} and \eqref{mix2} gives that
\begin{align}
\| u\|_{L^2 (B_{r_1})}
&\le C  |\log r_1| r_1 \brac{\frac{\pr{K+|\log r_0|} F_\de\pr{r_0}}{r_0} \|u\|_{L^2(B_{2r_0})}}^{k_0} \nonumber \\ 
&\times \brac{\frac{\pr{K+|\log R_1|} F_\de\pr{R_1}}{R_1} \|u\|_{L^2(B_{R_1})}}^{1 - k_0} \nonumber \\
&+C \pr{\frac{R_1}{r_0}} \pr{1 +\frac{|\log r_0|}{K}}e^{\pr{1 + C_1 K^\kappa + C_2 M^\mu}
\pr{\phi\pr{\frac{R_1}{2}}-\phi(r_0)}} \|u\|_{L^2(B_{2r_0})}.
 \label{end2}
\end{align}
By elliptic regularity (see for example \cite{HL11}, \cite{GT01}) and a scaling argument, we have that
\begin{align}
\|u\|_{L^\infty(B_r)}
\le C \frac{F_\de\pr{r}}{r} \|u\|_{L^2(B_{2r})}.
\label{ell}
\end{align}
From \eqref{end2} and \eqref{ell}, we arrive at the three-ball inequality in the $L^\infty$-norm that is given in \eqref{three}.
\end{proof}

Now we are ready to give the proof of Theorem \ref{thh}.
We first use the three-ball inequality to perform the propagation of smallness argument.
Then we apply the three-ball inequality again to obtain the order of vanishing estimate.

\begin{proof} [Proof of Theorem  \ref{thh}]
Without loss of generality, we may assume that $x_0 = 0$.
Let $\eps > 0$ be sufficiently small.
Fix some $\de > 0$.
Choose $r_0=\frac{r}{2}$, $r_1=4r$ and $R_1=10r$.
The application of \eqref{three} shows that
\begin{align}
\|u\|_{L^\infty \pr{B_{3r}}}
&\le C F_\de\pr{1}^2 \pr{K+|\log r|} |\log r|  \|u\|_{L^\iny(B_{r})}^{k_0} \|u\|_{L^\iny(B_{10r})}^{1 - k_0} \nonumber \\
&+C \pr{1 + \frac{\abs{\log r}}{K}} \exp\brac{ \pr{1 + C_1 K^\kappa + C_2 M^\mu} \pr{\phi\pr{5r}-\phi\pr{\frac r 2}} } \|u\|_{L^\iny(B_{r})},
\label{refi}
\end{align}
where $\disp k_0 = \frac{\phi(5r)-\phi(4r)}{\phi(5r)-\phi\pr{\frac r
2}}$. It is easy to check that
$$c^{-1}\leq \phi(5r)-\phi\pr{\frac{r}{2}}\leq c \quad \mbox{and} \quad c^{-1}\leq \phi(5r)-\phi(4r)\leq c,$$
where $c$ is some universal constant.
Thus, $k_0$ is independent of $r$ in this case.

Choose a small $r$ such that
$$\|u \|_{L^\iny\pr{B_r}}=\ell,$$
where $\ell>0$ by the unique continuation property.
Since $\|u \|_{L^\iny\pr{B_1}} \geq 1$, there exists some $\bar x\in B_1$ such that $\disp u(\bar x)=\|u \|_{L^\iny\pr{B_1}}\geq 1$.
We select a sequence of balls with radius $r$, centered at $x_0=0, \ x_1, \ldots, x_d$ so that $x_{i+1}\in B_{r}(x_i)$ for every $i$, and $\bar x\in B_{r}(x_d)$.
Notice that the number of balls, $d$, depends on the radius $r$ which will be fixed later.
Applying the $L^\infty$ version of three-ball inequality (\ref{refi}) at the origin and using boundedness assumption of $u$ given in \eqref{bound} yields that
\begin{align*}
\|u\|_{L^\infty \pr{B_{3r}(0)}}
&\le C \hat C^{1 - k_0} F_\de\pr{1}^2  \ell^{k_0}  \pr{1+\frac{|\log r|}{K}} \abs{\log r} K  \nonumber \\
&+ C \ell \pr{1 + \frac{|\log r|}{K}} \exp\brac{c\pr{1 + C_1 K^\kappa + C_2 M^\mu} }.
\end{align*}
Since $B_r(x_{i+1})\subset B_{3r}(x_{i})$, it is true that
\begin{equation}
\|u\|_{L^\infty (B_r(x_{i+1}))}\leq  \|u\|_{L^\infty
(B_{3r}(x_{i}))} \label{bbb}
\end{equation}
for every $i = 1, 2, \ldots, d$.
Repeating the argument as before with balls centered at $x_i$ and using \eqref{bbb}, we get
\begin{equation*}
\|u\|_{L^\infty (B_{3r}(x_{i}))}
\leq C_i \ell^{D_i} \pr{1+\frac{|\log r|}{K}}^{E_i} |\log r|^{F_i}  \exp\brac{H_i\pr{1 + C_1 K^\kappa + C_2 M^\mu} }
\end{equation*}
for $i=0, 1, \cdots, d$, where each $C_i$ depends on $d$, $\hat C$, $s$, $t$, $K$, $M$, and $C$ from Lemma \ref{threeBall}, and $D_i$, $E_i$, $F_i$ $H_i$ are constants that depend on $d$. Due to the fact that $u(\bar x)\geq 1$ and $\bar x \in B_{3r}(x_d)$, we have that
\begin{equation*}
\ell
\ge c \exp\brac{-C \pr{1 + C_1 K^\kappa + C_2 M^\mu} }\pr{1+\frac{|\log r|}{K}}^{-C} |\log r|^{-C},
\end{equation*}
where $c$ and $C$ are new constants with $c\pr{s, t, d, K, M, \hat C, \eps}$ and $C\pr{d}$.

Now we fix the radius $r$ as a small number.
In this sense, $d$ is a fixed constant.
We are going to use the three-ball inequality at the origin again with a different set of radii.
Let $\frac{3}{4}r_1=r$, $R_1=10r$ and let $r_0 << r$, i.e. $r_0$ is sufficiently small with
respect to $r$.
It follows from the three-ball inequality \eqref{three} that,
$$ \ell \leq {\rm I} +\Pi,$$
where
\begin{align*}
{\rm I} &= C F_\de\pr{r} |\log r| \brac{ (K+|\log r_0|) F_\de\pr{r_0}\|u\|_{L^\iny(B_{2r_0})}}^{k_0}
 \brac{ (K+|\log 10r|) F_\de\pr{10 r} \|u\|_{L^\iny(B_{10 r})}}^{1 - k_0} \\
\Pi &= C F_\de\pr{r} \pr{\frac{r }{r_0 }} \pr{1 +\frac{|\log r_0|}{K}} e^{\pr{1 + C_1 K^\kappa + C_2 M^\mu} \pr{\phi\pr{5r}-\phi(r_0)}} \|u\|_{L^\iny(B_{2r_0})},
\end{align*}
with $\disp k_0 = \frac{\phi(5r)-\phi(\frac 4 3 r)}{\phi(5r)-\phi(r_0)}$ and $F_\de\pr{r} = 1 + r K^{\frac{s}{s-2} + \de} + r M^{\frac{t}{2t-2}+\de}$.

On one hand, if ${\rm I} \leq \Pi$, we have
\begin{align*}
&c \exp\brac{-C \pr{1 + C_1 K^\kappa + C_2 M^\mu} }\pr{1+\frac{|\log r|}{K}}^{-C} |\log r|^{-C}
\le \ell \le 2 \Pi \\
&\le 2 C F_\de\pr{r} \pr{\frac{r }{r_0 }} \pr{1 +\frac{|\log r_0|}{K}} e^{\pr{1 + C_1 K^\kappa + C_2 M^\mu} \pr{\phi\pr{5r}-\phi(r_0)}} \|u\|_{L^\iny(B_{2r_0})}.
\end{align*}
If $r_0 << r$, then $\phi\pr{r_0}-\pr{C + \phi\pr{5r}} \gtrsim
\phi\pr{r_0}$. Since $r$ is a fixed small positive constant, we
obtain
\begin{align*}
\|u\|_{L^\iny(B_{2r_0})}
&\ge c r_0^{C \pr{1 + C_1 K^\kappa + C_2 M^\mu} },
\end{align*}
where now $c\pr{s, t, M, K, \hat C, \eps}$ and $C$ is some universal constant.

On the other hand, if $\Pi \leq {\rm I}$, then
\begin{align*}
&c \exp\brac{-C \pr{1 + C_1 K^\kappa + C_2 M^\mu} }\pr{1+\frac{|\log r|}{K}}^{-C} |\log r|^{-C}
\le \ell \le 2 {\rm I} \\
&\le 2 C F_\de\pr{r} |\log r| \brac{(K+|\log r_0|) F_\de\pr{r_0} \|u\|_{L^\iny(B_{2r_0})}}^{k_0} \brac{(K+|\log 10r|)  F_\de\pr{10 r} \|u\|_{L^\iny(B_{10 r})}}^{1 - k_0}.
\end{align*}
Raising both sides to $\frac{1}{k_0}$ and using that
$\|u\|_{L^\infty\pr{B_{10r}}}\leq \hat{C}$ gives that
\begin{align*}
 \|u\|_{L^\iny(B_{2r_0})}
 &\ge \frac{\hat C}{|\log r_0|} \pr{\frac{c/2 C \hat C K}{\abs{\log r}^{1+C} \pr{K+|\log r|}^C F_\de\pr{1}^{2} }}^{\frac 1 {k_0}} \exp\brac{-\frac{C}{k_0} \pr{1 + C_1 K^\kappa + C_2 M^\mu} }.
\end{align*}
Since $r$ is a fixed small positive constant, then
$\frac{1}{k_0}\simeq\log \frac{1}{r_0}$ if $ r_0\ll r$. Finally, we
arrive at
\begin{align*}
 \|u\|_{L^\iny(B_{2r_0})}
 &\ge C r_0^{C\pr{1 + C_1 K^\kappa + C_2 M^\mu}}
\end{align*}
as before.
This completes the proof of Theorem \ref{thh}.
\end{proof}

\section{The proof of Theorem \ref{thhh}}
\label{Thm3Proof}

The proof of Theorem \ref{thhh} is in the same spirit as that of Theorem \ref{thh}.
The main difference between these theorems is that Theorem \ref{thhh} has smaller values for $\mu$ than those in Theorem \ref{thh}.
See the remark following the statement of Theorem \ref{thhh} for the specific comparison of powers.
The improvement in Theorem \ref{thhh} comes from a slightly modified Carleman estimate for the operator $\LP + V$ which is made possible by the absence of a gradient potential.

\begin{lemma}
Assume that for some $t \in \pb{ 1, \ \infty}$, $\norm{V}_{L^t\pr{B_{R_0}}} \le M$.
Then for every sufficiently small $\eps > 0$, there exist constants $C_0$, $C_1$, and
sufficiently small $R_0 < 1$  such that for any $u\in C^{\infty}_{0}(B_{R_0}(x_0)\setminus \set{x_0})$ and
$$\tau \ge 1+ C_1 M^{\mu},$$
one has
\begin{align}
\tau^{2 - \frac 1 p} \|(\log r)^{-1} e^{-\tau \phi(r)}u\|_{L^2(r^{-2}dx)}
&\leq  C_0 \|(\log r ) e^{-\tau \phi(r)} r^2\pr{ \LP u + V u}\|_{L^p(r^{-2} dx)} ,
\label{main12}
\end{align}
where
$\disp \mu = \left\{\begin{array}{ll}
\frac{2t}{3t-2} & 2 < t \le \iny \medskip \\
\frac{t}{2t-2-\eps\pr{2t-1-2\eps}} & 1 < t \le 2
\end{array}\right.$,
and
$\disp p = \left\{\begin{array}{ll}
\frac{2t}{t+2} & t > 2 \medskip \\
\frac{t}{t- \eps} & 1 < t \le 2
\end{array}\right.$.
Moreover,
$C_0 = 2C$, where $C\pr{p, q, \eps}$ is from Theorem \ref{Carlpq}
with $ q = \frac{t}{t-1-\eps}$ when $1 < t \le 2$, and $C_1=
C_1\pr{t, \eps}$. \label{CarlpqV}
\end{lemma}

\begin{proof}
Assume that $\eps$ is sufficiently small.
The same argument as \eqref{triIneq} from the proof of Lemma \ref{CarlpqVW} implies that
\begin{align}
&\tau^{2 - \frac 1 p} \|(\log r)^{-1} e^{-\tau \phi(r)}u\|_{L^2(r^{-2}dx)}
+\tau^{\frac 2 q \pr{1 - \eps} + 1 - \frac 1 p} \|(\log r)^{-1} e^{-\tau \phi(r)}u\|_{L^q(r^{-2\pr{1 - \eps}}dx)} \nonumber \\
&\le C\|(\log r) e^{-\tau \phi(r)} r^2 (\LP u + V u)\|_{L^p(r^{-2} dx)}
+ C\|(\log r) e^{-\tau \phi(r)} r^2 V u \|_{L^p(r^{-2} dx)}.
\label{triIneqV}
\end{align}
We proceed to discuss the cases for $t$ to obtain (\ref{main12}).
\\

\nid {\bf Case 1: $t \in \pb{2, \iny}$} \\
Set $p = \frac{2t}{t+2}$ so that $t=\frac{2p}{2-p}$.
Substituting the estimate \eqref{hod3} in \eqref{triIneqV} yields that
\begin{align*}
&\tau^{\frac 3 2 - \frac 1 t} \|(\log r)^{-1} e^{-\tau \phi(r)}u\|_{L^2(r^{-2}dx)} \\
&\le C\|(\log r) e^{-\tau \phi(r)} r^2 (\LP u + V u)\|_{L^p(r^{-2} dx)}
+ c C M \|(\log r)^{-1} e^{-\tau \phi(r)} u\|_{L^2(r^{-2}dx)}.
\end{align*}
Taking $\tau \ge 1 + \pr{2 c C M}^{\frac{2t}{3t-2}}$, we are able to absorb the second term on the righthand side of the last inequality into the lefthand side.
Then the estimate (\ref{main12}) follows.
\\

\nid {\bf Case 2:  $t \in \pb{1, 2}$} \\
We optimize the power of $\tau$ by choosing $p$ very close to $1$.
Set $p =\frac{1}{1-\frac{\eps}{t}}$ and $q=\frac{pt}{t-p}$.
Since $t \le 2$, then $q=\frac{t}{t-1-\eps}>2$.
And if $\eps < t -1$, $q < \iny$ as well.
We will use the estimate \eqref{hod33} to bound the last term in \eqref{triIneqV}.
We have
$$\frac{2}{q}(1-\eps)+1-\frac{1}{p}=\frac{1-2\eps}{p}-\frac{2(1-\eps)}{t}+1=\frac{2t-2-\eps\pr{2t-1-2\eps}}{t}>0$$
if we further assume that $\eps\pr{2t-1-2\eps} < 2\pr{t-1}$.
Substituting \eqref{hod33} into \eqref{triIneqV} and simplifying shows that
\begin{align*}
&\tau^{1 + \frac{\eps}{t}} \|(\log r)^{-1} e^{-\tau \phi(r)}u\|_{L^2(r^{-2}dx)}
+\tau^{\frac{ 2t-2-\eps\pr{2t-1-2\eps}}{t}} \|(\log r)^{-1} e^{-\tau \phi(r)}u\|_{L^q(r^{-2\pr{1 - \eps}}dx)} \nonumber \\
&\le C\|(\log r) e^{-\tau \phi(r)} r^2 (\LP u + Vu)\|_{L^p(r^{-2} dx)}
+ c C M \|(\log r)^{-1} e^{-\tau \phi(r)} u\|_{L^q(r^{-2\pr{1 - \eps}}dx)}.
\end{align*}
We may absorb the second term on the right if $\tau \ge 1 +\pr{c C M}^{\frac t {2t-2 -\eps\pr{2t-1-2\eps}}}$.
This completes the proof of Lemma \ref{CarlpqV}.
\end{proof}

\begin{proof}[Proof of Theorem  \ref{thhh}]
Using Lemma \ref{CarlpqV} in place of Lemma \ref{CarlpqVW}, we derive a $L^\infty$ version of the three-ball inequality like the one in Lemma \ref{threeBall} for the operator $\LP + V$.
Then by applying the propagation of smallness argument from the proof of Theorem \ref{thh}, we will arrive at the proof of Theorem \ref{thhh}.
\end{proof}

\section{Unique continuation at infinity}
\label{QuantUC}

In this section, using the scaling arguments established in \cite{BK05}, we show how the maximal order of vanishing estimates implies the quantitative unique continuation estimates at infinity.

\begin{proof}[Proof of Theorem \ref{UCVW}]
Let $u$ be a solution to \eqref{goal} in $\R^2$.
Fix $x_0 \in \R^2$ and set $\abs{x_0} = R$.
We do a scaling as follows, $u_R(x) = u(x_0 + Rx)$, $W_R\pr{x} = R \, W\pr{x_0 + R x}$, and $V_R\pr{x} = R^2 V\pr{x_0 + R x}$.
For any $r > 0$,
\begin{align*}
\norm{W_R}_{L^s\pr{B_r\pr{0}}}
&= \pr{\int_{B_r\pr{0}} \abs{W_R\pr{x}}^s dx}^{\frac 1 s}
= \pr{\int_{B_r\pr{0}} \abs{R \, W\pr{x_0 + R x}}^s dx}^{\frac 1 s} \\
&= R^{1 - \frac 2 s } \pr{  \int_{B_r\pr{0}} \abs{W\pr{x_0 + R x}}^s
d\pr{Rx} }^{\frac 1 s} = R^{1 - \frac 2 s } \norm{W}_{L^s\pr{B_{r
R}\pr{x_0}}}.
\end{align*}
and
\begin{align*}
\norm{V_R}_{L^t\pr{B_r\pr{0}}}
&= R^{2 - \frac 2 t } \pr{  \int_{B_r\pr{0}} \abs{V\pr{x_0 + R x}}^t
d\pr{Rx} }^{\frac 1 t} = R^{2 - \frac 2 t } \norm{V}_{L^t\pr{B_{r
R}\pr{x_0}}}.
\end{align*}
Therefore,
$$\norm{W_R}_{L^s\pr{B_{10}\pr{0}}} = R^{1 - \frac 2 s} \norm{W}_{L^s\pr{B_{10R}\pr{x_0}}} \le A_1 R^{1 - \frac 2 s}$$
and
$$\disp \norm{V_R}_{L^t\pr{B_{10}\pr{0}}} \le A_0 R^{2 - \frac 2 t}.$$
Moreover, $u_R$ satisfies a scaled version of \eqref{goal} in $B_{10}$,
\begin{align*}
& \LP u_R\pr{x} + W_R\pr{x} \cdot \gr u_R\pr{x}  + V_R\pr{x} u_R\pr{x} \\
&= R^2 \LP u\pr{x_0 + R x} + R^2 W\pr{x_0 + Rx} \cdot \gr u\pr{x_0 + R x}  + R^2 V\pr{x_0 + R x}u\pr{x_0 + R x}
= 0.
\end{align*}
Clearly,
\begin{align*}
\norm{u_R}_{L^\iny\pr{B_{6}}}
&= \norm{u}_{L^\iny\pr{B_{6R}\pr{x_0}}} \le C_0.
\end{align*}
Note that for $\disp\widetilde{x_0} := -x_0/R$, we have $\disp \abs{\widetilde{x_0}} = 1$ and  $\abs{u_R(\widetilde{x_0})} = \abs{u(0)} \ge 1$.
Thus, $\disp\norm{u_R}_{L^\iny(B_1)} \ge 1$.
If $R$ is sufficiently large, then we may apply Theorem \ref{thh} to $u_R$ with some arbitrarily small $\eps \in \pr{ 0, 1}$, $K = A_1 R^{1 - \frac 2 s}$, $M = A_2 R^{2 - \frac 2 t}$, and $\hat C = C_0$ to obtain
\begin{align*}
\norm{u}_{L^\iny\pr{{B_{1}(x_0)}}} = & \norm{u_R}_{L^\iny\pr{B_{1/R}(0)}}  \\
\ge & c(1/R)^{^{C\brac{1 + C_1 \pr{A_1 R^{1 - \frac 2 s}}^\kappa + C_2 \pr{A_2 R^{2 - \frac 2 t} }^\mu}}} \\
=& c \exp\set{-C\brac{1 + C_1 \pr{A_1 R^{1 - \frac 2 s}}^\kappa +
C_2 \pr{A_2 R^{2 - \frac 2 t} }^\mu} \log R},
\end{align*}
where $\kappa$ and $\mu$ depend on $s$, $t$, and $\eps$.
Further simplifying, we see that
\begin{align*}
\norm{u}_{L^\iny\pr{{B_{1}(x_0)}}}
\ge &c \exp\brac{-C\pr{1 + C_1 A_1^\kappa + C_2 A_2^\mu } R^\Pi \log R}
\end{align*}
where
\begin{align*}
\Pi := \max\set{\kappa \pr{\frac {s-2} s}, \mu \pr{\frac {2t-2} t }}.
\end{align*}
Recalling the values of  $\kappa$ and $\mu$ from Theorem \ref{thh}, a computation shows that when $0 < \eps < 1$,
\begin{align*}
\Pi =
\left\{\begin{array}{ll}
2 & t > \frac{2s}{s+2}, \medskip \\
 \frac{t(s-2)}{s(t-1-\eps t)} & 1 < t \le \frac{2s}{s+2} \\
\end{array}\right.
\end{align*}
and the conclusion of the theorem follows.
\end{proof}

\begin{proof}[Proof of Theorem \ref{UCV}]
 To prove Theorem \ref{UCV}, we follow the same approach as before and
use that
$$\mu \pr{\frac {2t-2} t}  =  \left\{\begin{array}{ll}
\frac{4t-4}{3t-2} & t > 2 \medskip \\
\frac{2t-2}{2t-2-\eps\pr{2t-1-2\eps}} & 1 < t \le 2,
\end{array}\right.$$
where $\mu$ is given in the statement of Theorem \ref{thhh}.
\end{proof}

\renewcommand{\theequation}{A.\arabic{equation}}

\setcounter{equation}{0}  
\section*{Appendix}  

In this section, we first present the proof of Lemma \ref{CarLpp} for the case $1<p<2$.
The eigenfunction estimates in Lemma \ref{upDown} play an important role in the argument.
Then we prove a quantitative Caccioppoli inequality in dimension $n=2$.

\begin{proof}[Proof of Lemma \ref{CarLpp}]
 We define a conjugated operator $L^-_\tau$ of $L^-$ by
$$L^-_\tau u=e^{-\tau\varphi(t)}L^-(e^{\tau \varphi(t)}u).$$
With $v=e^{\tau \varphi(t)}u$, it is equivalent to prove
\begin{equation}
\|t^{-{1}}u \|_{L^2(dtd\omega)}\leq C \tau^\beta \|t L^-_\tau
u\|_{L^p(dtd\omega)}.
\label{keyu}
\end{equation}
From the definition of $\Lambda$ and $L^-$ in \eqref{ord} and \eqref{use}, the operator $L^-_\tau$ can be written as
\begin{equation}
L^-_\tau=\sum_{k\geq 0} (\partial_t+\tau \varphi'(t)-k)P_k.
\label{ord1}
\end{equation}
Let $M=\lceil 2\tau\rceil$.
Since $\disp \sum_{k\geq 0} P_k v= v$, we split $\disp \sum_{k \ge 0} P_k v$ into
$$P^+_\tau=\sum_{k> M}P_k, \quad \quad  P^-_\tau=\sum_{k=0}^{M}P_k.$$
Then \eqref{keyu} is reduced to show both the following inequalities
\begin{equation}
\| t^{-1} P^+_\tau {u} \|_{L^2(dtd\omega)} \leq \tau^{\beta}\|
t{L_\tau^- u} \|_{L^{{p}}(dtd\omega)} \label{key1}
\end{equation}
and
\begin{equation}
\| t^{-1} P^-_\tau {u}\|_{L^2(dtd\omega)}\leq \tau^\beta\| t
{L_\tau^- u} \|_{L^{{p}}(dtd\omega)} \label{key2}
\end{equation}
hold for all $u \in C^\iny_c\pr{(-\infty, \ t_0)\times S^{1}}$ and $1< p < 2$.
We first establish \eqref{key1}.
From \eqref{ord1} and properties of the projection operator $P_k$, it follows that
\begin{equation}
P_k L^-_\tau u= (\partial_t+\tau \varphi'(t)-k)P_k u. \label{sord}
\end{equation}
For $u\in C^\infty_{0}\pr{ (-\infty, \ t_0)\times S^{1}}$, the solution $P_k u$ of this first order differential equation is given by
\begin{equation}
P_k u(t, \omega)
=-\int_{-\infty}^{\infty} H(s-t)e^{k(t-s)+\tau\brac{\varphi(s)-\varphi(t)}} P_k L^-_\tau u (s,\omega)\, ds, \label{star}
\end{equation}
where $H(z)=1$ if $z\geq 0$ and $H(z)=0$ if $z<0$.

For $k\geq M \ge 2 \tau$, it can be shown that
\begin{equation*}
H(s-t)e^{k(t-s)+\tau\brac{\varphi(s)-\varphi(t)}} \leq
e^{-\frac{1}{2}k|t-s|}
\end{equation*}
for all $s, t\in (-\infty, \ t_0)$.
Taking the $L^2\pr{S^{1}}$-norm in \eqref{star} yields that
\begin{equation*}
\|P_k u(t, \cdot)\|_{L^2(S^{1})}
\leq  \int_{-\infty}^{\infty} e^{-\frac{1}{2}k|t-s|} \|P_k L^-_\tau u(s, \cdot)\|_{L^2(S^{1})}
\,ds.
\end{equation*}
Applying the eigenfunction estimates \eqref{indu} gives that
\begin{equation*}
\|P_k u(t, \cdot)\|_{L^2(S^{1})}
\leq C \int_{-\infty}^{\infty} e^{-\frac{1}{2}k|t-s|} \|L^-_\tau u(s, \cdot)\|_{L^p(S^{1})} \,ds
\end{equation*}
for all $1\leq p\leq 2$.
Furthermore, the application of Young's inequality for convolution yields that
\begin{equation*}
\|P_k u\|_{L^2(dt d\omega)}
\leq C  \pr{\int_{-\infty}^{\infty} e^{-\frac{\sigma}{2}k|z|} dz}^{\frac{1}{\sigma}}\|L^-_\tau u\|_{L^p(dtd\omega)}
\end{equation*}
with $\frac{1}{\sigma}=\frac{3}{2}-\frac{1}{p}$.
Therefore,
\begin{equation*}
\|P_k u\|_{L^2(dt d\omega)}
\leq C k^{\frac{1}{p} - \frac{3}{2}} \|L^-_\tau u\|_{L^p(dtd\omega)},
\end{equation*}
where we have used the fact that
$$ \pr{\int_{-\infty}^{\infty} e^{-\frac{\sigma}{2}k|z|} dz }^{\frac{1}{\sigma}}\leq C k^{\frac{1}{p} - \frac{3}{2}}.  $$
Squaring and summing up $k> M$ shows that
$$\sum_{k> M} \|P_k u\|^2_{L^2(dt d\omega)} \leq C \sum_{k> M} k^{\frac{2}{p} - 3} \|L^-_\tau
u\|^2_{L^p(dtd\omega)}.$$
Thus, $\disp \sum_{k> M} k^{\frac{2}{p} - 3}$ converges if $p>1$.
If $p=1$, $\disp \sum_{k> M} k^{\frac{2}{p}- 3}$ diverges.
Therefore,
\begin{equation*}
\|  P^+_\tau u\|_{L^2(dtd\omega)}
\leq C\tau^{\frac{1-p}{p}}\| L^-_\tau u\|_{L^p(dtd\omega)},
\end{equation*}
which implies estimate \eqref{key1} since $|t|\geq |t_0|$, where $\abs{t_0}$ is large.

Set $N=\lceil\tau \varphi^\prime(t)\rceil$.
Recall that $\varphi(t)=t+\log t^2$.
By Taylor's theorem, for all $s, t \in (-\infty, \ t_0)$, we have
\begin{align}
\varphi(s)-\varphi(t) & = \varphi'(t)(s-t)-\frac{1}{(s_0)^2}(s-t)^2,
\label{Taylor}
\end{align}
where $s_0$ is some number between $s$ and $t$.
If $s>t$, then
$$S_k(s, t) = e^{k(t-s) + \tau\brac{\vp\pr{s - \vp\pr{t}}}} \leq e^{-(k-\tau\varphi'(t))(s-t)-\frac{\tau}{t^2}(s-t)^2}.$$
Hence
\begin{equation}
H(s-t) S_k(s, t)\leq e^{-|k- N 
||s-t|-\frac{\tau}{t^2}(s-t)^2}.
\label{case1}
\end{equation}
Furthermore, we consider the case $N\leq k\leq M$.
The summation of \eqref{star} over $k$ shows that
\begin{equation}
\| \sum^M_{k=N} P_k u(t, \cdot)\|_{L^2(S^{1})} \leq
\int_{-\infty}^{\infty} \| \sum^M_{k=N} H(s-t)S_k(s,t) P_k L^-_\tau u(s, \cdot) \|_{L^2(S^{1})}\, ds. \label{back}
\end{equation}
Let $c_k= H(s-t)S_k(s,t)$.
It is clear that $|c_k|\leq 1$.
Now we make use of Lemma \ref{upDown}.
An application of estimate \eqref{haha} shows that for all $1< p < 2$
\begin{align}
& \|\sum^M_{k=N}H(s-t)S_k(s,t) P_k L^-_\tau u(s, \cdot)\|_{L^2(S^{1})}\leq \nonumber \medskip\\
& C \pr{\sum^M_{k=N} H(s-t)|S_k(s,t)|^2}^{\frac{1}{p}-\frac{1}{2}}
 \|L^-_\tau u(s, \cdot) \|_{L^p(S^{1})}.
\label{mixSum}
\end{align}
From \eqref{case1}, we have
\begin{align}
\sum^M_{k=N} H(s-t)|S_k(s,t)|^2 &\leq \pr{\sum^M_{k=N+1} e^{-2|k-
N||s-t|}+1}e^{ -\frac{2\tau}{t^2}(s-t)^2}
\nonumber \\
& \leq  C \pr{ \frac{1}{|s-t|}+1} e^{ -\frac{\tau}{t^2}(s-t)^2} .
\label{sumBnd}
\end{align}
 Therefore, the last two inequalities imply that
\begin{align*}
\|\sum^M_{k=N}H(s-t)S_k(s,t) P_k L^-_\tau u(s,\cdot)\|_{L^2(S^{1})} 
& \leq  C (|s-t|^{-\alpha_2}+1)e^{-\frac{\alpha_2\tau}{t^2}(s-t)^2}\| L^-_\tau u(s, \cdot) \|_{L^p(S^{1})},
\end{align*}
where  $\alpha_2=\frac{(2-p)}{2p}$.
It can be shown that
\begin{equation}
e^{-\frac{\alpha_2\tau}{t^2}(s-t)^2}\leq
C|t|\pr{1+\sqrt{\tau}|s-t|}^{-1}. \label{expBnd}
\end{equation}
Thus,
$$ \|\sum^M_{k=N}H(s-t)S_k(s,t) P_k L^-_\tau u(s,\cdot)\|_{L^2(S^{1})}
\leq \frac{C|t|(|s-t|^{-\alpha_2}+1)\| L^-_\tau u(s, \cdot) \|_{L^p(S^{1})}}{1+\sqrt{\tau}|s-t|}.$$
It follows from \eqref{back} that
\begin{equation}
|t|^{-1}\| \sum^M_{k=N} P_k u(t, \cdot)\|_{L^2(S^{1})} \leq C
 \int_{-\infty}^{\infty}\frac{(|s-t|^{-\alpha_2}+1)\| L^-_\tau u(s, \cdot) \|_{L^p(S^{1})}}{1+\sqrt{\tau}|s-t|}.
\label{mile}
\end{equation}

For the case $k\leq N-1$, solving the first order differential equation \eqref{sord} gives that
\begin{equation}
P_k u(t, \omega)=\int_{-\infty}^{\infty} H(t-s)S_k(s, t) P_k
L^-_\tau u\pr{s, \om}\, ds \label{star1}.
\end{equation}
The estimate \eqref{Taylor} shows that for any $s, t$
\begin{equation}
H(t-s) S_k(s, t)\leq e^{-|N- 1 - k ||s-t|-\frac{\tau}{s^2}(t-s)^2}.
\label{case2}
\end{equation}
Using \eqref{case2} and performing the calculation as before, we conclude that
\begin{equation}
\| \sum_{k=0}^{N-1} P_k u(t, \cdot)\|_{L^2(S^{1})} \leq C
\int_{-\infty}^{\infty} \frac{|s|(|s-t|^{-\alpha_2}+1)\| L^-_\tau u(s, \cdot) \|_{L^p(S^{1})}}{1+\sqrt{\tau}|s-t|}.
\label{mile1}
\end{equation}
Since $s,t $ are in $(-\infty, \ t_0)$ with $|t_0|$ large enough, the combination of estimates \eqref{mile} and \eqref{mile1} gives
\begin{equation*}
|t|^{-1}\| P_\tau^- u(t, \cdot)\|_{L^2(S^{1})} \leq C
 \int_{-\infty}^{\infty}
\frac{|s|(|s-t|^{-\alpha_2}+1)\| L^-_\tau u(s, \cdot)
\|_{L^p(S^{1})}}{1+\sqrt{\tau}|s-t|}.
\end{equation*}
Applying Young's inequality for convolution, we obtain
\begin{equation*}
\| t^{-1} P_\tau^- u (t, \cdot) \|_{L^2(dtd\omega)} \leq C
\brac{\int_{-\infty}^{\infty} \pr{\frac{ |z|^{-\alpha_2}+1}{1+\sqrt{\tau}|z|}}^\si \, dz}^{\frac{1}{\si}} \| tL^-_\tau u \|_{L^p(dtd\omega)},
\end{equation*}
where $\frac{1}{\si}=\frac{3}{2}-\frac{1}{p}$.
Therefore,
\begin{equation*}
\| t^{-1} P_\tau^- u (t, \cdot) \|_{L^2(dtd\omega)}\leq
C\tau^{-\frac{1}{2\si}+\frac{\alpha_2}{2}} \| t L^-_\tau u \|_{L^p(dtd\omega)},
\end{equation*}
where we have used the fact
$$\brac{\int_{-\infty}^{\infty} \pr{\frac{ |z|^{-\alpha_2}+1} {1+\sqrt{\tau}|z|}}^\si \, dz}^{\frac{1}{\si}}
\leq C\tau^{-\frac{1}{2\si}+\frac{\alpha_2}{2}}$$
with $\al_2 \in \pr{0, \frac 1 2}$ and $\si \in \pr{1, 2}$.
 This completes
\eqref{key2} since $-\frac{1}{2\si}+\frac{\alpha_2}{2}=\frac{1-p}{p}$.
Finally, the proof is complete.
\end{proof}

We state and prove a Caccioppoli inequality for the second order elliptic equation (\ref{goal}) with singular lower order terms.
Because our lower order terms are assumed to be singular, we must employ a Sobolev embedding in $\R^2$, and this forces the right hand side to be relatively larger than it was in Lemma 5 from \cite{DZ17}, the corresponding result for $n \ge 3$.

\begin{lemma}
Assume that for some $s \in \pb{2, \iny}$ and $t \in \pb{ 1, \iny}$, $\norm{W}_{L^s\pr{B_{R}}} \le K$ and $\norm{V}_{L^t\pr{B_{R}}} \le M$.
Let $u$ be a solution to equation \eqref{goal} in $B_R$.
Then for any $\de > 0$, there exists a constant $C$, depending only on $s$, $t$ and $\de$, such that for any $r<R$,
\begin{equation*}
\|\nabla u\|^2_{L^2(B_r)}\leq
C\brac{\frac{1}{(R-r)^2}+M^{\frac{t}{t-1}+ \de} +K^{\frac{2s}{s-2} + \de }}\| u\|^2_{L^2(B_R)}.
\end{equation*}
\label{CaccLem}
\end{lemma}

\begin{proof}
We start by decomposing $V$ and $W$ into bounded and unbounded parts.
For some $M_0, K_0$ to be determined, let
$$V(x)=\overline{V}_{M_0}+V_{M_0}, \quad W(x)=\overline{W}_{K_0}+W_{K_0},$$
where
$${\overline{V}}_{M_0}=V(x)\chi_{\set{|V(x)|\leq \sqrt{M_0}}}, \quad
{V}_{M_0}=V(x)\chi_{\set{|V(x)|>\sqrt{M_0}}},$$ and
$${\overline{W}}_{K_0}=W(x)\chi_{\set{|W(x)|\leq \sqrt{K_0}}}, \quad
{W}_{K_0}=W(x)\chi_{\set{|W(x)|>\sqrt{K_0}}}.$$

For any $q \in \brac{1, t}$,
\begin{equation}
\|V_{M_0}\|_{L^q}\leq
M_0^{-\frac{t-q}{2q}}\|V_{M_0}\|_{L^{t}}^{\frac{t}{q}}\leq
M_0^{-\frac{t-q}{2q}}\|V\|_{L^{t}}^{\frac{t}{q}} \leq
M_0^{-\frac{t-q}{2q}}M^{\frac{t}{q}}. \label{who}
\end{equation}
Similarly, for any $q \in \brac{1, s}$, we have
\begin{equation}
\|W_{K_0}\|_{L^q}\leq
K_0^{-\frac{s-q}{2q}}\|W_{K_0}\|_{L^{s}}^{\frac{s}{q}}\leq
K_0^{-\frac{s-q}{2q}}\|W\|_{L^{s}}^{\frac{s}{q}} \leq
K_0^{-\frac{s-q}{2q}}K^{\frac{s}{q}}. \label{mmo}
\end{equation}

Let $B_R\subset B_1$.
Choose a smooth cut-off function $\eta \in C^\iny_0\pr{B_R}$ such that $\eta(x) \equiv 1$ in $B_r$ and $|\nabla \eta|\leq \frac{C}{|R-r|}$.
Multiplying both sides of equation (\ref{goal}) by $\eta^2 u$ and integrating by parts, we obtain
\begin{equation}
\int |\nabla u|^2 \eta^2 = \int V\eta^2 u^2 +\int W\cdot \nabla u \,
\eta^2 u - 2\int \nabla u\cdot\nabla \eta \, \eta \, u. \label{cacc}
\end{equation}
We estimate the terms on the right side of \eqref{cacc}.
For the first term, we see that
\begin{equation*}
\int V \eta^2 u^2 \leq \int |{\overline{V}_{M_0}}|\eta^2 u^2 + \int
|{{V}_{M_0}}|\eta^2 u^2.
\end{equation*}
It is clear that
\begin{equation}
\int |{\overline{V}_{M_0}}|\eta^2 u^2\leq M_0^{\frac{1}{2}}\int
\eta^2 u^2. \label{sss}
\end{equation}
Fix $\de > 0$.
Let $\de_0 = \de \frac{\pr{t-1}^2}{t+\de\pr{t-1}}$.
By H\"older's inequality, \eqref{who} with $q = 1 + \de_0$, and Sobolev embedding with $2q^\prime := 2\pr{1 + \de_0^{-1}} > 2$, we get
\begin{align}
\abs{\int {V}_{M_0}\eta^2 u^2 } &\leq
\pr{\int|{{V}_{M_0}}|^{q}}^{\frac1{q}} \pr{\int
\abs{\eta^2u^2}^{q^\prime}}^{\frac{1}{q^\prime}} \le C_{\de_0}
{M_0}^{-\frac{t-q}{2q}}M^{\frac{t}{q}} \int |\nabla (\eta u)|^2.
\label{mmm}
\end{align}
Taking $C_{\de_0} {M_0}^{-\frac{t-q}{2q}}M^{\frac{t}{q}}=\frac 1 {16}$, i.e $M_0 = C_{\de_0, t} M^{\frac{2t}{t-q}}$, from (\ref{sss}) and (\ref{mmm}), we get
\begin{eqnarray}
\int V \eta^2 u^2 & \leq & C M^{\frac{t}{t-q}}\int |\eta
u|^2+\frac{1}{16}\int |\nabla (\eta u)|^2 \nonumber
\medskip \\
& \leq &C M^{\frac{t}{t-1} + \de}\int |\eta u|^2+\frac{1}{8}\int
|\nabla \eta|^2 u^2 +\frac{1}{8}\int |\nabla u|^2 \eta^2.
\label{more1}
\end{eqnarray}

Now we estimate the second term in the righthand side of (\ref{cacc}).
We have that
\begin{equation}
\int W \cdot \nabla u \eta^2 u = \int |{\overline{W}_{K_0}}| |\nabla
u| \eta^2 u + \int |{W}_{K_0}| | \nabla u| \eta^2 u. \label{w1}
\end{equation}
By Young's inequality for products,
\begin{equation}
\int |{\overline{W}_{K_0}}|  |\nabla u| \eta^2 u \leq \frac{1}{8}
\int |\nabla u|^2\eta^2 +CK_0 \int \eta^2 u^2. \label{w2}
\end{equation}
This time, set $\de_0 =\de \frac{\pr{s-2}^2}{2s+\de\pr{s-2}}$.
By H\"older's inequality, \eqref{who} with $q = 2 + \de_0$, and Sobolev embedding with $q^\prime = 2\pr{1 + 2 \de_0^{-1}} > 2$, we get
\begin{eqnarray}
\int {W}_{K_0} \cdot \nabla u \eta^2 u
&\leq & \pr{\int |{{W}_{K_0}}|^{{q}}}^{\frac{1}{q}} \pr{\int |\nabla u\cdot \eta |^{2}}^{\frac{1}{2}} \pr{\int |u \eta |^{q^{\prime}}}^{\frac{1}{q^\prime}} \nonumber \medskip \\
& \leq &C_{\de_0}  K_0^{-\frac{s-q}{2q}}K^{\frac{s}{q}}\||\nabla u|\eta\|_{L^2} \|\nabla(\eta u)\|_{L^2} \nonumber \medskip \\
& \leq &C_{\de_0}  K_0^{-\frac{s-q}{2q}}K^{\frac{s}{q}}\||\nabla u|\eta\|_{L^2}\pr{\||\nabla\eta | u \|_{L^2} + \||\nabla u|\eta\|_{L^2} }\nonumber \medskip \\
&\leq & 2 C_{\de_0} K_0^{-\frac{s-q}{2q}}K^{\frac{s}{2}}\||\nabla
u|\eta|\|^2_{L^2} + \frac 1 2 C_{\de_0}
K_0^{-\frac{s-q}{2q}}K^{\frac{s}{2}} \||\nabla\eta | u \|_{L^2}^2.
\label{w3}
\end{eqnarray}
We choose $C_{\de_0} K_0^{-\frac{s-q}{2q}}K^{\frac{s}{q}}=\frac 1 {16}$, that is, $K_0=C_{\de_0, s} K^{\frac{2s}{s-q}}$.
The combination of (\ref{w1}), (\ref{w2}) and (\ref{w3}) gives that
\begin{eqnarray}
\int W \cdot \nabla u \, \eta^2 u &\leq& \frac{1}{4} \||\nabla u|
\eta\|^2_{L^2}+CK^{\frac{2s}{s-2} + \de}\|u\eta\|^2_{L^2}
+\frac{1}{32}\||\nabla \eta| u\|^2_{L^2} . \label{more2}
\end{eqnarray}
Finally, Young's inequality for products implies that
\begin{align*}
2\int \nabla u\cdot\nabla \eta \, \eta \, u \le \frac 1 8  \||\nabla u| \eta\|^2_{L^2} + 8  \||\nabla \eta| u \|^2_{L^2}
\end{align*}

Together with (\ref{cacc}), (\ref{more1}) and (\ref{more2}), we obtain
\begin{equation*}
\int |\nabla u|^2 \eta^2 \leq C \pr{M^{\frac{t}{t-1} + \de}
+K^{\frac{2s}{s-2} + \de}} \int |u\eta|^2\,dx+ C \int |\nabla
\eta|^2u^2\,dx.
\end{equation*}
From the assumptions on $\eta$, this completes the proof in the lemma.
\end{proof}


\def\cprime{$'$}

\end{document}